 \documentclass[final,3p,times]{elsarticle}

\usepackage{amssymb}
 \usepackage{amsthm}
\usepackage{amsmath,amssymb,amsopn,amsfonts,mathrsfs,amsbsy,amscd}

\newcommand{\prs}{\langle\;,\;\rangle}
\newcommand{\too}{\longrightarrow}
\newcommand{\om}{\omega}
\newcommand{\esp}{\quad\mbox{and}\quad}

\newcommand{\G}{{\mathfrak{g}}}

\newcommand{\ad}{{\mathrm{ad}}}

\newcommand{\Om}{\Omega}
\newcommand{\na}{\nabla}

\newcommand{\al}{\alpha}
\newcommand{\be}{\beta}
\newcommand{\ga}{\gamma}

\newcommand{\e}{\epsilon}

\newcommand{\De}{\Delta}
\newcommand{\de}{\delta}

\font\bb=msbm10

\def\K{\hbox{\bb K}}
\def\Z{\hbox{\bb Z}}

\def\R{\hbox{\bb R}}

\newtheorem{Def}{Definition}[section]
\newtheorem{theo}{Theorem}[section]
\newtheorem{pr}{Proposition}[section]

\newtheorem{co}{Corollary}[section]

\newtheorem{exem}{Example}
\newtheorem{remark}{Remark}

\begin{document}

\begin{frontmatter}


 \title{Left invariant para-K\"ahler and hyper-para-K\"ahler structures on Lie groups\tnoteref{label1}}
 
 \fntext[label3]{This research was conducted within the framework of Action concert\'ee CNRST-CNRS Project SPM04/13.}

\title{On para-K\"ahler and hyper-para-K\"ahler Lie algebras}

 \author[label1,label2]{Sa\"id Benayadi, Mohamed Boucetta}
 \address[label1]{Universit\'e de Lorraine, Laboratoire IECL, CNRS-UMR 7502,\\ Ile du Saulcy, F-57045 Metz
 cedex
 1, France.\\e-mail: said.benayadi@univ-lorraine.fr
 }
 \address[label2]{Universit\'e Cadi-Ayyad\\
 Facult\'e des sciences et techniques\\
 BP 549 Marrakech Maroc\\e-mail: boucetta@fstg-marrakech.ac.ma
 }



\begin{abstract}We study  Lie algebras admitting para-K\"ahler and 
hyper-para-K\"ahler structures. We give new characterizations of these Lie algebras
and we develop many methods to build large classes of examples.  Bai considered para-K\"ahler Lie algebras as left symmetric
bialgebras. We reconsider this point of view and improve it in order to obtain some new results. The study of para-K\"ahler and 
hyper-para-K\"ahler is intimately linked to the study of left symmetric algebras and, in particular, those admitting invariant symplectic forms. In this paper, we give many new classes of left symmetric algebras and a complete description of all associative algebras admitting an invariant symplectic form. We give also all four dimensional hyper-para-K\"ahler Lie algebras.

\end{abstract}


\end{frontmatter}






\section{Introduction}\label{section1}
A  \emph{para-complex} structure on a $2n$-dimensional manifold $M$
is a field $K$ of involutive endomorphisms $(K^2 = Id_{TM})$ such that 
the eigendistributions $T^\pm M$ with eigenvalues $\pm1$ have constant rank $n$ and are
integrable.
In the presence of
a pseudo-Riemannian metric this notion leads to the notion of para-K\"ahler manifolds. A
\emph{para-K\"ahler} structure on a manifold $M$ is a pair $(g,K)$ where $g$ is a
pseudo-Riemannian metric and $K$ is a parallel skew-symmetric para-complex
structure. If $(g,K)$ is a para-K\"ahler structure on $M$, then $\om = g \circ
K$ is a symplectic
structure and the $\pm1$-eigendistributions $T^\pm M$ of $K$ are two integrable
$\om$-Lagrangian distributions. Due to this, a para-K\"ahler structure can be
identified with a bi-Lagrangian structure $(\om, T^\pm M)$ where $\om$ is a symplectic
structure and $T^\pm M$ are two integrable Lagrangian distributions. If $(M,g,K)$ is 
a para-K\"ahler manifold and $J$ is a parallel field  of skew-symmetric endomorphisms
such that $J^2=-Id_{TM}$ and $JK=-KJ$ then $(M,g,K,J)$ is called a
\emph{hyper-para-K\"ahler} manifold or \emph{hyper-symplectic} manifold.
The notion of almost para-complex
structure (or almost product structure) on a manifold was introduced
by P.K. Rasevskii \cite{ra} and P. Libermann \cite{liberman}. The
paper \cite{cru} contains a survey on  para-K\"ahler geometries.
Hyper-para-K\"ahler structures were introduced by N. Hitchin in \cite{hitchin} and have
become an important subject of study
lately, due mainly to their applications in theoretical physics (specially in dimension 4).
See
for instance \cite{pope}, where there is a discussion on the relationship between
hyper-para-K\"ahler metrics and the $N = 2$ superstring. When the manifold is a Lie group
$G$,  the metric and the para-complex structure are considered left-invariant, they
are both determined by their restrictions to the Lie algebra $\G$ of $G$.  In a such
situation, $(\G,g_e,K_e)$ is called   \emph{para-K\"ahler  Lie
algebra}. We recover also the notion of \emph{hyper-para-K\"ahler Lie algebra} when we
start from  a left invariant hyper-para-K\"ahler structure on a Lie group. Para-K\"ahler and hyper-para-K\"ahler Lie algebras has been studied by many authors \cite{andrada, salomon, bai, bai1, bajo}.\\
This paper is
devoted to the study of para-K\"ahler and hyper-para-K\"ahler Lie algebras. We present some known results by adopting a new approach which we think simplify both the presentation and the proofs. In the large part of the paper, we give some new results which permit a better understanding of this algebras and
  the construction of a rich classes
of non trivial new examples.
  The basic tools of the study of
para-K\"ahler and hyper-para-K\"ahler Lie algebras are two
types of algebras: left symmetric algebras which have been studied by many authors and a
less known class, namely, left symmetric algebras endowed with  invariant symplectic forms called in \cite{bai0} special symplectic Lie algebras.
We call such algebras symplectic left symmetric algebras. Our study leads incidentally to the construction of a large classes of left symmetric algebras,  symplectic left symmetric algebras and symplectic Lie algebras.\\

 The paper is organized as follows. In Sections \ref{section2} and \ref{section2bis} we recall some basic
definitions and we present, by using a new approach, some known characterizations of para-K\"ahler Lie algebras. We give a particular attention to
  the
notion of left symmetric bialgebras
 introduced by Bai \cite{bai}. We introduce the notion of quasi $S$-matrices as a generalization of $S$-matrices introduced by Bai. Proposition \ref{pralg} describing the Lie algebra structure of the para-K\"ahler Lie algebra associated to a quasi $S$-matrix  will play a crucial role in Sections \ref{section5}-\ref{section5bis}. It shows also (see Remark \ref{remdiatta} $(b)$) that a quasi $S$-matrix on a left symmetric algebra $U$ defines a Lie triple system on $U^*$ (see \cite{jacobson, lister, smirnov} for the definition and properties of Lie triple systems).
In Section
\ref{section4}, we develop some general methods to build new examples of para-K\"ahler Lie algebras. 
In Section \ref{section3}, we give a new characterization of
 hyper-para-K\"ahler Lie algebras based on  a notion of compatibility 
 between two left symmetric products on a given vector space (see  Theorem \ref{theohyperkahler} and Definition \ref{Def1}).
   Sections \ref{section5}-\ref{section5bis} are devoted to the study of quasi $S$-matrices on symplectic Lie algebras,  on symplectic left symmetric algebras and on pseudo-Riemannian flat Lie algebras.
On a symplectic Lie algebra with its canonical left symmetric product the set of quasi $S$-matrices is in bijection with the set of solutions of an equation which generalizes the modified Yang-Baxter equation (see Proposition \ref{prqss}). As a consequence, we find a  method to build a new  class of para-K\"ahler Lie algebras (see Theorem \ref{theosymp}) and actually a new class of Lie algebras with a Lie triple system (see Remark \ref{rem3}).   On a symplectic left symmetric algebra or a pseudo-Riemannian flat Lie algebra, the set of $S$-matrices is in bijection with the set of operators generalizing $\mathcal{O}$-operators (see Proposition \ref{prqsss}). As a consequence, we find a  method to build a new  class of para-K\"ahler and hyper-para-K\"ahler Lie algebras (see Theorems \ref{theohyper}-\ref{theohyper2}). We get also a new class of Lie algebras with a Lie triple system (see Remark \ref{rem4}).
 In Section \ref{section6} we give all
four dimensional para-K\"ahler Lie algebras. We use a method which is different from the
one used in \cite{andrada} and which has the advantage of simplifying
enormously the calculations. We devote Section \ref{section7} to  a complete description  of associative
symplectic left symmetric algebras (see Theorems \ref{main1}-\ref{main2}).

\paragraph{Notations} For a Lie algebra $\G$, its bracket will be denoted by $[\;,\;]$ and for any $u\in\G$, $\ad_u$ is the endomorphism of $\G$ given by $\ad_u(v)=[u,v]$. If $A:\G\too\G$ is an endomorphism, the Nijenhuis torsion of $A$ is given by
\begin{equation}\label{nui}N_A(u,v):=[Au,Av]-A[Au,v]-A[u,Av]+A^2[u,v].
\end{equation}

If $(U,.)$ is an algebra, for any $u\in U$, $\mathrm{L}_u,\mathrm{R}_u:U\too U$ denote the left and the right multiplication by $u$ given by $\mathrm{L}_u(v)=u.v$ and $\mathrm{R}_u(v)=v.u$.  The commutator of $(U,.)$ is the bracket on $U$ given by $[u,v]=u.v-v.u$. The curvature of $(U,.)$ is the tensor $\mathrm{K}$ given by 
\begin{equation*}\label{eq12b}
 \mathrm{K}(u,v):=[\mathrm{L}_{u},\mathrm{L}_{v}]-\mathrm{L}_{[u,v]}.
\end{equation*}
Then, for any $u,v,w\in U$, we have the Bianchi identity 
 \begin{equation}\label{bianchi}\oint[u,[v,w]]=\oint \mathrm{K}(u,v)w,
 \end{equation}where $\oint$ stands of the cyclic sum.
The product on $U$ is called Lie-admissible if its commutator is a Lie bracket, i.e., for any $u,v,w\in U$,
\begin{equation*}\oint[u,[v,w]]=\oint \mathrm{K}(u,v)w=0.
 \end{equation*}
Let  $V$ be a vector space and $F$ is an endomorphism of $V$. We denote by  $F^t:V^*\too V^*$  the dual endomorphism. For any $X\in V$ and $\al\in V^*$ we denote $\al(X)$ by $\prec\al,X\succ$.
The phase space of $V$ is the vector space $\Phi(V):=V\oplus V^*$ endowed with the two nondegenerate bilinear forms $\prs_0$ and $\Om_0$ given by
\begin{equation*}\label{fdual}
\langle u+\al,v+\be\rangle_0=\prec\al,v\succ+\prec\be,u\succ\esp
\Om_0(u+\al,v+\be)=\prec\be,u\succ-\prec\al,v\succ.
\end{equation*}We denote by $K_0:\Phi(V)\too\Phi(V)$ the endomorphism given by $ K_0(u+\al)=u-\al$.  \\
Let $\om\in\wedge^2V^*$ which is nondegenerate. We denote by $\flat:V\too V^*$ the isomorphism given by $\flat(v)=\om(v,.)$. Put $T(V):=V\times V$ and define $\prs_1,\Om_1,K_1,J_1$ on $T(V)$ by
\begin{equation*} \Om_1[(u,v),(w,z)]=\om(z,u)-\om(v,w),\; \langle(u,v),(w,z)\rangle_1=\om(z,u)+\om(v,w),\end{equation*}
\begin{equation*} K_1(u,v)=(u,-v)
\esp J_1(u,v)=(-v,u). \end{equation*}
Finally, if $\rho:\G\too\mathrm{End}(V)$ is a representation of a Lie algebra, we denote by ${\rho}^*:\G\too\mathrm{End}(V^*)$ the dual representation given by $\rho^*(X)(\al)=-\rho(X)^t(\al)$.

\section{Some definitions}\label{section2}
In this section, we recall the definitions of different types of algebraic structures  used through this paper.
 \begin{itemize}\item 
 
 A complex structure on a Lie 	algebra	$\G$ is an isomorphism $J:\G\too\G$ satisfying $J^2=-\mathrm{Id}_\G$ and $N_J=0$. A complex structure $J$ is called abelian if, for any $u,v\in\G$,
 \[ [Ju,Jv]=[u,v]. \]
 A para-complex structure on a Lie 	algebra	$\G$ is an isomorphism $K:\G\too\G$ satisfying $K^2=\mathrm{Id}_G$, $N_K=0$ and $\dim\ker(K+\mathrm{Id}_\G)=
 \dim\ker(K-\mathrm{Id}_\G)$. A para-complex structure $K$ is called abelian if, for any $u,v\in\G$,
  \[ [Ku,Kv]=-[u,v]. \]
  A complex product structure on $\G$ is a couple $(J,K)$ where $J$ is a complex structure, $K$ is a para-complex structure and $KJ=-JK$.

\item  A \emph{pseudo-Riemannian Lie algebra} is a finite dimensional Lie algebra 
$(\G,[\;,\;])$ 
endowed with a bilinear symmetric nondegenerate form $\prs$. The associated 
\emph{Levi-Civita
product} is the product on $\G$, $(u,v)\mapsto u.v$, given by  Koszul's formula
\begin{equation}\label{eq11}2\langle
u.v,w\rangle=\langle[u,v],w\rangle+\langle[w,u],v\rangle+\langle[w,v],u\rangle.
\end{equation}
This product is entirely determined by the facts that it is Lie-admissible, i.e.,
$[u,v]=u.v-v.u$ and,
for any $u\in\G$, the left multiplication by $u$ is skew-symmetric with respect to $\prs$. We call $(\G,[\;,\;])$ flat when its Levi-Civita product has a vanishing curvature.
\item An algebra $(U,.)$ is called \emph{left symmetric } if
\begin{equation*}\label{eq12} \mathrm{ass}(u,v,w)=\mathrm{ass}(v,u,w),\end{equation*}where
$\mathrm{ass}(u,v,w)=(u.v).w-u.(v.w)$. This relation is equivalent to the vanishing of
the curvature of $(U,.)$. The relation \eqref{bianchi} implies that a left symmetric product is Lie-admissible.
 If $(U,.)$ is a left symmetric algebra then the Lie algebra $(U,[\;,\;])$ has
two representations, namely, $\ad_U:U\too\mathrm{End}(U),$ $u\mapsto\ad_u$ and 
$\mathrm{L}_U:U\too\mathrm{End}(U),$ $u\mapsto \mathrm{L}_{u}$.
\\
Any associative algebra is a left symmetric algebra and   if a left symmetric algebra is
abelian then it is associative.
\item A \emph{symplectic Lie algebra} is a Lie algebra $(\G,\om)$ endowed
with a bilinear skew-symmetric nondegenerate form $\om$ such that, for any $u,v,w\in\G$,
\begin{equation*}\label{eq14}
 \om([u,v],w)+\om([v,w],u)+\om([w,u],v)=0.
\end{equation*}
It is a well-known result \cite{chu} that the product ${\mathbf{a}}:\G\times\G\too\G$ given by
\begin{equation}\label{eq15}
 \om({\mathbf{a}}(u,v),w)=-\om(v,[u,w])
\end{equation}induces on $\G$ a left symmetric algebra structure satisfying
${\mathbf{a}}(u,v)-{\mathbf{a}}(v,u)=[u,v]$. We call ${\mathbf{a}}$ the
\emph{left symmetric product} associated to $(\G,\om)$.
\item  A \emph{symplectic left symmetric algebra} is a left symmetric algebra
$(U,.)$ endowed
with a 
bilinear skew-symmetric nondegenerate form $\om$ which is invariant, i.e., for any $u,v,w\in U$,
$$\om(u.v,w)+\om(v,u.w)=0.$$
This implies that $(U,[\;,\;],\om)$ is a symplectic Lie algebra.
Symplectic left symmetric algebras, called in \cite{bai1} special symplectic Lie algebras,  play a central role in the study of
hyper-para-K\"ahler Lie algebras (see Section \ref{section3}). 
\end{itemize}
\section{Para-K\"ahler Lie algebras}\label{section2bis}

 The notion of para-K\"ahler Lie algebra is very subtle  and has many equivalent definitions depending on if one emphasizes on its pseudo-Riemannian metric and the associated Levi-Civita product or on its symplectic form and its associated left symmetric product. It has also many characterizations. There is at least three such characterizations in \cite{bai}: a para-K\"ahler Lie algebra can be characterized as the phase space of a Lie algebra, as a matched pair of Lie algebras or as a left symmetric bialgebra. In this section, we choose 	 the pseudo-Riemannian point of view and we give a new characterization based on Bianchi identity  \eqref{bianchi}. This characterization has the advantage of leading easily to the notions of left symmetric bialgebra and $S$-matrix  introduced in \cite{bai}. To end this section, we introduce a generalization of the notion of $S$-matrix and we give a precise description of the  para-K\"ahler Lie algebras associated to these generalized $S$-matrices.
 
 A \emph{para-K\"ahler Lie algebra} is a pseudo-Riemannian Lie algebra $(\G,\prs)$ endowed with
 an isomorphism $K:\G\too \G$ satisfying
 $K^2=\mathrm{Id}_\G$, $K$  is skew-symmetric with respect to $\prs$ and
$K$ is invariant with respect to the Levi-Civita product, i.e., 
$\mathrm{L}_u\circ K=K\circ \mathrm{L}_u$, for any $u\in\G$. A para-K\"ahler Lie
algebra $(\G,\prs,K)$  carries a natural  bilinear skew-symmetric nondegenerate form
$\Om$ defined by
$\Om(u,v)=\langle Ku,v\rangle,$ and one can see easily that:
\begin{enumerate}\item $(\G,K)$ is a para-complex Lie algebra,
\item $(\G,\Om)$ is a symplectic Lie algebra,
 \item $\G=\G^1\oplus\G^{-1}$ where $\G^1=\ker(K-\mathrm{Id}_{\G})$ and
$\G^{-1}=\ker(K+\mathrm{Id}_{\G})$,
\item $\G^1$ and $\G^{-1}$ are subalgebras isotropic with respect to $\prs$ and 
Lagrangian with 
respect to $\Om$,
 \item for
any $u\in\G$, $u.\G^1\subset\G^1$ and $u.\G^{-1}\subset\G^{-1}$ (the dot is the
Levi-Civita product).
\end{enumerate}
A para-K\"ahler Lie algebra  carries two products: the Levi-Civita
product and the left symmetric product $\mathbf{a}$ associated to $(\G,\Om)$. The following proposition
clarifies the relations between them. Its proof is an easy computation. 

\begin{pr}\label{pr11}Let $(\G,\prs,\Om,K)$ be a para-K\"ahler Lie algebra. Then,
   for any $u,v\in \G^{1}$ and $\al,\be\in \G^{-1}$,
 $$u.v=\mathbf{a}(u,v)\quad\mbox{and}\quad\al.\be=\mathbf{a}(\al,\be).$$
 In particular, $\G^1$ and $\G^{-1}$ are left symmetric algebras.
\end{pr}

Let $(\G,\prs,\Om,K)$ be a para-K\"ahler Lie algebra. For any $u\in\G^{-1}$, let
$u^*$ denote the element of $(\G^{1})^*$ given by
$\prec u^*,v\succ=\langle u,v\rangle.$ The map $u\mapsto u^*$ realizes an isomorphism between 
$\G^{-1}$ and $(\G^{1})^*$. Thus we can identify $(\G,\prs,\Om,K)$ to the phase space $(\Phi(\G^1),\prs_0,\Om_0,K_0)$.  With this identification in mind  
and according to Proposition \ref{pr11},  the Levi-Civita product induces a product
on $\G^1$ and 
$(\G^{1})^*$ which coincides with the affine product $\mathbf{a}$. Thus both $\G^1$ and
$(\G^{1})^*$ carry a
left symmetric algebra structure.  For any $u\in\G^1$ and for any $\al\in(\G^{1})^*$,
we denote by $\mathrm{L}_u:\G^1\too\G^1$ and $\mathrm{L}_\al:(\G^{1})^*\too(\G^{1})^*$
the left multiplication by $u$ and $\al$, respectively, i.e., for any $v\in\G^1$ and
any $\be\in(\G^{1})^*$,
$$\mathrm{L}_uv=u.v=\mathbf{a}(u,v)\esp
\mathrm{L}_\al\be=\al.\be=\mathbf{a}(\al,\be).$$The right
multiplications $\mathrm{R}_u$ and $\mathrm{R}_\al$ are defined in a similar manner.
The following proposition shows that
the Levi-Civita product  and the affine product  on $\G$ identified with
$\Phi(\G^1)$ are entirely determined by their restrictions to $\G^1$ and
$(\G^{1})^*$. The proof of this proposition is straightforward.
\begin{pr}\label{pr12}Let $\G$ be a para-K\"ahler Lie algebra identified with 
$\Phi(\G^1)$ as above. Then
\begin{enumerate}\item For any $u\in \G^1$ and for any $\al\in (\G^{1})^*$,
$$u.\al=-\mathrm{L}_u^t\al\quad\mbox{and}\quad 
\al. u=-\mathrm{L}_\al^tu.$$
\item For any $u\in \G^{1}$ and for any $\al\in (\G^{1})^*$,
$$\mathbf{a}(u,\al)=\mathrm{R}_\al^tu-\ad_u^t\al
\quad\mbox{
and}
\quad
\mathbf{a}(\al, u)=-\ad_\al^tu+\mathrm{R}_u^t\al,$$where $\ad_u:\G^1\too\G^1$ and
$\ad_\al:(\G^{1})^*\too(\G^{1})^*$ are given by $\ad_uv=[u,v]$ and
$\ad_\al\be=[\al,\be]$.

\end{enumerate}

\end{pr}

Conversely, let $U$ be a finite dimensional vector space and $U^*$ its dual space. We
suppose that both $U$
and $U^*$ have a 
structure of left symmetric algebra.  We extend the products on $U$ and $U^*$ to
$\Phi(U)$
by putting,
 for any $X,Y\in U$ and for any $\al,\be\in U^*$,
\begin{equation}\label{eq16}(X+\al).(Y+\be)=X.Y
-\mathrm{L}_\al^tY-\mathrm{L}_X^t\be+\al.\be.\end{equation} 
We consider the two bilinear maps
$\rho:U\times U^*\too\mathrm{End}(U)$ {and}
$\rho^*:U^*\times U\too\mathrm{End}(U^*)$ defined by
\begin{equation}\label{eq17}\rho(X,\al)=[\mathrm{L}_X,\mathrm{L}_\al^t]
+\mathrm{L}_{\mathrm{L}^t_\al X}+\mathrm{L}_{\mathrm{L}^t_X\al}^t\esp\rho^*(\al,X)=[\mathrm{L}_\al,\mathrm{L}_X^t]
+\mathrm{L}_{\mathrm{L}^t_X \al}+\mathrm{L}_{\mathrm{L}^t_\al X}^t.\end{equation}
Note that the endomorphism $\rho^*(\al,X)$ is the dual of $\rho(X,\al)$. We give now a new characterization of para-K\"ahler Lie algebras using Bianchi identity.
\begin{pr}\label{pr13}With the hypothesis above,
   the product on $\Phi(U)$ given   by \eqref{eq16} is Lie-admissible if and
only if 
  \begin{equation}\label{eq19}\rho(X,\al)Y=\rho(Y,\al)X\quad\mbox{and}\quad
\rho^*(\al,X)\be=\rho^*(\be,X)\al\end{equation}
for any
  $X,Y\in U$ and any $\al,\be\in U^*$. Moreover, this product  is left symmetric if and
  only if 
    $\rho(X,\al)=0,$ for any
    $X\in U$ and any $\al\in U^*$. 

\end{pr}

{\bf Proof.} According to Bianchi identity \eqref{bianchi}, the product given by \eqref{eq16} is
Lie-admissible if and only if, for any $u,v,w\in\Phi(U)$,
$$\mathrm{K}(u,v)w+\mathrm{K}(v,w)u+\mathrm{K}(w,u)v=0,$$
where $\mathrm{K}$ is the curvature of the product.
Since the products on $U$ and $U^*$ are left symmetric, for any $X,Y,Z\in U$ and for any $\al,\be,\ga\in U^*$,
$$\mathrm{K}(X,Y)Z=\mathrm{K}(\al,\be)\ga=0.$$ Thus the Bianchi identity is equivalent to
\begin{eqnarray*}
 \mathrm{K}(X,Y)\al+\mathrm{K}(Y,\al)X+\mathrm{K}(\al,X)Y=0,\quad\quad(*)\\
 \mathrm{K}(X,\al)\be+\mathrm{K}(\al,\be)X+\mathrm{K}(\be,X)\al=0,\quad\quad(**)
\end{eqnarray*}for any $X,Y\in U$ and for any $\al,\be\in U^*$. Now, we have obviously,
$$\mathrm{K}(X,Y)\al=\left([\mathrm{L}_X,\mathrm{L}_Y]-\mathrm{L}_{[X,Y]}\right)^t\al
=0.$$On the other hand, a direct computation gives
$$\mathrm{K}(Y,\al)X=\rho(Y,\al)X\esp \mathrm{K}(\al,X)Y=-\rho(X,\al)Y.$$
So $(*)$ is equivalent to $\rho(X,\al)Y=\rho(Y,\al)X$. A similar computation shows that
$(**)$ is equivalent to $\rho^*(\al,X)\be=\rho^*(\be,X)\al$. The second part of the
proposition follows easily from what above.\hfill$\square$\\

\begin{Def}
Two left symmetric products on  $U$ and $U^*$ satisfying \eqref{eq19} will be
called \emph{Lie-extendible}.\end{Def}
Thus we get the following result.
\begin{theo}\label{theoextendible}
Let $(U,.)$ and $(U^*,.)$ two Lie-extendible left symmetric products. Then $(\Phi(U),\prs_0,K_0)$ endowed with the Lie algebra bracket associated to the  product given by \eqref{eq16} is a para-K\"ahler Lie algebra. Moreover, all para-K\"ahler Lie algebras are obtained in this way.

\end{theo}
\begin{exem}\label{exemple1}
Let $(U,.)$ be a left symmetric algebra. Then the left symmetric product on $U$ and the trivial left symmetric product on $U^*$ are Lie-extendible so $(\Phi(U),\prs_0,K_0)$ endowed with the Lie algebra bracket associated to the left symmetric product
\begin{equation}\label{exem1eq1} (X+\al)\triangleright(Y+\be)=X.Y-\mathrm{L}_X^t\be \end{equation}is a para-K\"ahler Lie algebra. We denote by $[\;,\;]^\triangleright$ the Lie bracket associated to $\triangleright$. We have 
\[ [X+\al,Y+\be]^\triangleright=[X,Y]-\mathrm{L}_X^t\be+\mathrm{L}_Y^t\al. \]
This is just the semi-direct product of $(U,[\;,\;])$ with $U^*$ endowed with the trivial bracket and the action of $U$ on $U^*$ is given by $\mathrm{L}_U^*$. Moreover, it is easy to check that $(\Phi(U),[\;,\;]^\triangleright,\prs_0)$ is a flat pseudo-Riemannian  Lie algebra
and $(\Phi(U),\triangleright,\Om_0)$ is a  symplectic left symmetric algebra (see also Proposition 4.3 \cite{bai1}).
\end{exem}

In \cite{bai}, Bai gave a characterization of para-K\"ahler Lie algebras similar to the one used for Lie
bialgebras. He called these structures  \emph{left symmetric bialgebras}.  Let us present this point of view
in a new way by using Proposition \ref{pr13}.\\
We consider two left symmetric algebras $(U,.)$ and $(U^*,.)$. The  products on $U$ and $U^*$ define by duality, respectively, two maps $\mu:U^*\too
U^*\otimes
U^*$ and $\xi:U\too U\otimes U$. As Lie algebras $U$ and $U^*$ have two representations
 $\Psi_U:U\too \mathrm{End}(U\otimes U)$ and
$\Psi_{U^*}:U^*\too
\mathrm{End}(U^*\otimes U^*)$ given by
\begin{eqnarray*}
 \Psi_U=\mathrm{L}_U\otimes\ad_U\quad\mbox{and}\quad
\Psi_{U^*}=
\mathrm{L}_{U^*}\otimes\ad_{U^*}.
\end{eqnarray*}
A direct computation gives, for any $X,Y\in U$ and for any $\al,\be\in U^*$, 
\begin{eqnarray*}
\prec \be,\rho(X,\al)Y-\rho(Y,\al)X\succ&=&\Psi_U(X)(\xi(Y))(\al,\be)-\Psi_U(Y)(\xi(X))(\al,
\be)-\xi(
[ X , Y ]
)(\al,\be),\\
\prec\rho^*(\al,X)\be-\rho^*(\be,X)\al,Y\succ&=&\Psi_{U^*}(\al)(\mu(\be))(X,
Y)-\Psi_{U^*}(\be)(\mu(\al))(X , Y)-\mu( [
\al,\be]
)(X,Y).
\end{eqnarray*}By using Proposition \ref{pr13}, we recover a result of Bai (see \cite{bai} Theorem 4.1).
\begin{pr}\label{pr14}
 The product on $\Phi(U)$ given  by \eqref{eq16} is Lie-admissible if and only if $\xi$ is a
1-cocycle of
$(U,[\;,\;])$ with respect to the representation $\Psi_{U}$ and $\mu$ is a 1-cocycle of
$(U^*,[\;,\;])$ with respect to the representation $\Psi_{U^*}$, i.e., for any $X,Y\in U, \al,\be\in U^*$,
\begin{eqnarray*}
 \xi([X,Y])&=&\Psi_{U}(X)(\xi(Y))-\Psi_{U}(Y)(\xi(X)),\\
\mu([\al,\be])&=&\Psi_{U^*}(\al)(\mu(\be))-\Psi_{U^*}(\be)(\mu(\al)).
\end{eqnarray*}
\end{pr}
We  consider now the case where
 $\xi$ is a co-boundary. Indeed, 
let  $(U,.)$ be a left symmetric
algebra  and $\xi:U\too U\otimes U$ a
co-boundary of
$(U,[\;\;])$ with respect to $\Psi_{U}$, i.e., $\xi=\de\mathrm{r}$ where
 $\mathrm{r}\in U\otimes U$. By duality, $\xi$ define a product on $U^*$ by
\begin{equation}\label{eq110}\prec\al.\be,X\succ=\mathrm{r}(\mathrm{L}_X^t\al,\be)+\mathrm{r}
(\al,
\ad_X^t\be)=\mathrm{L}_X\mathrm{r
}(\al,\be)-\prec\mathrm{L}_{\mathrm{r}_\#(\al)}^t\be,X\succ,\end{equation}
where $\mathrm{r}_\#:U^*\too U$ is given by
$\prec\be,\mathrm{r}_\#(\al)\succ=\mathrm{r}(\al,\be)$ and
$\mathrm{L}_X\mathrm{r}(\al,\be)=\mathrm{r}(\mathrm{L}_X^t\al,\be)+\mathrm{r}(\al,\mathrm
{L}_X^t\be).$\\
According to Proposition \ref{pr13}, to get a para-K\"ahler Lie algebra structure on $\Phi(U)$, $(U^*,.)$ must be a left
symmetric algebra and the couple  $(U,.),$ $(U^*,.)$ must be Lie-extendible. Note that, since 
$\xi=\de\mathrm{r}$, the first equation in \eqref{eq19} holds. Let us  find out under which conditions   the second equation in \eqref{eq19} holds and $(U^*,.)$ is a left symmetric algebra.
Put $\mathrm{r}=\mathfrak{a}+\mathfrak{s}$ where $\mathfrak{a}$ is skew-symmetric and
$\mathfrak{s}$ is
symmetric and define $\mathrm{L}(\mathfrak{a})\in U^*\otimes U\otimes U$  by
$$\mathrm{L}(\mathfrak{a})(X,\al,\be)=\mathrm{L}_X\mathfrak{a}(\al,\be).$$
It follows immediately from \eqref{eq110} that, for any $\al,\be\in U^*$ and $X\in U$,
\begin{equation}\label{eq110bis}
\mathrm{L}_\al^tX=\mathrm{r}_\#\circ \mathrm{L}_X^t\al+[X,\mathrm{r}_\#(\al)]
\esp\prec[\al,\be],X\succ=\prec\mathrm{L}_{\mathrm{r}_\#(\be)}^t\al
-\mathrm{L}_{\mathrm{r}_\#(\al)}^t\be,X\succ+2\mathrm{L}_X\mathfrak{a}(\al,\be).
\end{equation}
We consider the two representations $\mathrm{Q}:U\too \mathrm{End}(U\otimes U\otimes U)$ 
and $\mathrm{P}:U\too \mathrm{End}(U^*\otimes U\otimes U)$ given by
\begin{eqnarray*}\mathrm{Q}&=&\mathrm{L}_{U}\otimes \mathrm{L}_{U}\otimes \ad_{U}\esp
\mathrm{P}=\mathrm{L}^*_{U}\otimes \mathrm{L}_{U}\otimes \mathrm{L}_{U}.\end{eqnarray*}
We define also $\De(\mathrm{r})\in U\otimes U\otimes U\simeq\mathrm{End}(U^*\otimes U^*,U)$ by
\begin{equation}\label{eq111}\De(\mathrm{r})(\al,\be)=
\mathrm{r}_\#([\al,\be])-[\mathrm{r}_\#(\al),\mathrm{r}
_\#
(\be)].\end{equation}
 \begin{pr}\label{prnew}For any $X,Y\in U$ and $\al,\be,\ga\in U^*$, we have
\begin{eqnarray*}
\prec \rho^*(\al,X)\be-\rho^*(\be,X)\al,Y\succ&=&2
P(X)(\mathrm{L}(\mathfrak{a}))(Y,\al,\be),\\
\mathrm{ass}(\al,\be,\ga)-\mathrm{ass}(\be,\al,\ga)&=&\prec\ga,-
\mathrm{Q}(X)(\De(\mathrm{r}))(\al,\be)+\mathrm{r}_\#\left(
\rho^*(\al,X)\be-\rho^*(\be,X)\al\right)\succ.
\end{eqnarray*}\end{pr}
{\bf Proof.} Let us compute first the associator of $\al,\be,\ga\in U^*$ with respect to the product given by \eqref{eq110}. We have, for any $X\in U$,
\begin{eqnarray*}
\prec\mathrm{ass}(\al,\be,\ga),X\succ&=&\prec\al.(\be.\ga),X\succ-\prec(\al.\be).\ga,X\succ\\
&=&\mathrm{r}(\mathrm{L}_X^t\al,\be.\ga)+\mathrm{r}(\al,\ad_X^t(\be.\ga))
-\mathrm{r}(\mathrm{L}_X^t(\al.\be),\ga)-\mathrm{r}(\al.\be,\ad_X^t\ga)\\
&=&\prec\be.\ga,\mathrm{r}_\#(\mathrm{L}_X^t\al)\succ+\prec\ad_X^t(\be.\ga),
\mathrm{r}_\#(\al)\succ-
\prec\ga,\mathrm{r}_\#(\mathrm{L}_X^t(\al.\be))\succ-\prec\ad_X^t\ga,
\mathrm{r}_\#(\al.\be)\succ\\
&=&\mathrm{r}(\mathrm{L}_{\mathrm{r}_\#(\mathrm{L}_X^t\al)}^t\be,\ga)+\mathrm{r}(\be,\ad_{\mathrm{r}_\#(\mathrm{L}_X^t\al)}^t\ga)+
\mathrm{r}(\mathrm{L}_{[X,\mathrm{r}_\#(\al)]}^t\be,\ga)+
\mathrm{r}(\be,\ad_{[X,\mathrm{r}_\#(\al)]}^t\ga)\\
&&-
\prec\ga,\mathrm{r}_\#(\mathrm{L}_X^t(\al.\be))\succ-\prec\ga,
[X,\mathrm{r}_\#(\al.\be)]\succ\\
&=&\prec\ga,\mathrm{r}_\#\left(\mathrm{L}_{\mathrm{r}_\#(\mathrm{L}_X^t\al)}^t\be\right)\succ+\prec\ga,[\mathrm{r}_\#(\mathrm{L}_X^t\al),\mathrm{r}_\#(\be)]\succ+
\prec\ga,\mathrm{r}_\#\left(\mathrm{L}_{[X,\mathrm{r}_\#(\al)]}^t\be\right)\succ\\&&+\prec\ga,[[X,\mathrm{r}_\#(\al)],\mathrm{r}_\#(\be)]\succ
-
\prec\ga,\mathrm{r}_\#(\mathrm{L}_X^t(\al.\be))\succ-\prec\ga,
[X,\mathrm{r}_\#(\al.\be)]\succ.
\end{eqnarray*}
On the other hand,
\begin{eqnarray*} Q(X)(\De(\mathrm{r}))(\al,\be)&=&[X,
\De(\mathrm{r})(\al,\be)]-\De(\mathrm{r})(\mathrm{L}_X^*\al,\be)-\De(\mathrm{r})(\al,\mathrm{L}_X^*\be)\\
&=&[X,\mathrm{r}_\#([\al,\be])]-[X,[\mathrm{r}_\#(\al),\mathrm{r}_\#
(\be)]]+\mathrm{r}_\#([\mathrm{L}_X^t\al,\be])-[\mathrm{r}_\#\left(\mathrm{L}_X^t\al\right),\mathrm{r}
_\#
(\be)]\\&&+\mathrm{r}_\#([\al,\mathrm{L}_X^t\be])-[\mathrm{r}_\#(\al),\mathrm{r}
_\#
\left(\mathrm{L}_X^t\be\right)].
\end{eqnarray*}So 
\[ \mathrm{ass}(\al,\be,\ga)-\mathrm{ass}(\be,\al,\ga)+
\prec\ga,Q(X)(\De(\mathrm{r}))(\al,\be)\succ=\prec\ga,\mathrm{r}_\#(A)\succ, \]where
\begin{eqnarray*}  A&=&\mathrm{L}_{\mathrm{r}_\#(\mathrm{L}_X^t\al)}^t\be-\mathrm{L}_{\mathrm{r}_\#(\mathrm{L}_X^t\be)}^t\al +\mathrm{L}_{[X,\mathrm{r}_\#(\al)]}^t\be
-\mathrm{L}_{[X,\mathrm{r}_\#(\be)]}^t\al
-\mathrm{L}_X^t([\al,\be])+[\mathrm{L}_X^t\al,\be]+[\al,\mathrm{L}_X^t\be].
\end{eqnarray*}By using the first relation in \eqref{eq110bis}, we get
\[ A=\\
=\mathrm{L}_{\mathrm{L}_\al^tX}^t\be
-\mathrm{L}_{\mathrm{L}_\be^tX}^t\al-\mathrm{L}_X^t([\al,\be])+[\mathrm{L}_X^t\al,\be]+[\al,\mathrm{L}_X^t\be]. \]
Now, according to \eqref{eq17}, we have
\begin{eqnarray*} \rho^*(\al,X)\be&=&[\mathrm{L}_\al,\mathrm{L}_X^t]\be
+\mathrm{L}_{\mathrm{L}^t_X \al}\be+\mathrm{L}_{\mathrm{L}^t_\al X}^t\be\\
&=&\al.(\mathrm{L}_X^t\be)-\mathrm{L}_X^t(\al.\be)
+({\mathrm{L}^t_X \al}).\be+\mathrm{L}_{\mathrm{L}^t_\al X}^t\be, \end{eqnarray*}
so $A=\rho^*(\al,X)\be-\rho^*(\be,X)\al$ and the second assertion follows.
On the other hand, by using the second relation in \eqref{eq110bis}, we get  for any $Y\in U$, 
\begin{eqnarray*}
\prec A,Y\succ&=&-2\mathrm{L}_{X.Y}\mathfrak{a}(\al,\be)+
2\mathrm{L}_{Y}\mathfrak{a}(\mathrm{L}_{X}^t\al,\be)
+2\mathrm{L}_{Y}\mathfrak{a}(\al,\mathrm{L}_{X}^t\be),
\end{eqnarray*}and the first assertion follows.\hfill$\square$

So we get the following result.
\begin{theo}\label{theo11}Let $(U,.)$ be a left symmetric algebra and $\mathrm{r}=\mathfrak{a}+\mathfrak{s}\in U\otimes U$.
 Then the product given by \eqref{eq110} is left
symmetric and the left symmetric products on $(U,U^*)$ are Lie-extendible 
if and only if
 for any $X\in U$ $$\mathrm{Q}(X)(\De(\mathrm{r}))=0\quad\mbox{and}\quad 
\mathrm{P}(X)(\mathrm{L}(\mathfrak{a}))=0.$$
In this case, $(\Phi(U),\prs_0,\Om_0,K_0)$ endowed with the Lie bracket associated to the product given by \eqref{eq16} is a para-K\"ahler Lie algebra.

\end{theo}
We have immediately the following corollary.
\begin{co}\label{co11}
 If  $\mathfrak{a}$ is $\mathrm{L}_U$-invariant, i.e. $\mathrm{L}(\mathfrak{a})=0$, 
then the product given by \eqref{eq110} is left
symmetric and  the left symmetric products on $(U,U^*)$ are Lie-extendible if and only if
 $\De(\mathrm{r})$ is $Q$-invariant.

\end{co}

Actually the statement of Theorem \ref{theo11} is the same as the one of Theorem 5.4 in \cite{bai}. To show this, let us investigate the relation between $\De(\mathrm{r})$ and $[[\mathrm{r},\mathrm{r}]]$ appearing in Bai's Theorem. Let $(U,.)$ be a left symmetric algebra and $\mathrm{r}=\sum_{i}a_i\otimes b_i$. In \cite{bai}, Bai defines $[[\mathrm{r},\mathrm{r}]]$ by
\[ [[\mathrm{r},\mathrm{r}]]=\mathrm{r}_{13}.\mathrm{r}_{12}
-\mathrm{r}_{23}.\mathrm{r}_{21}+[\mathrm{r}_{23},\mathrm{r}_{12}] 
-[\mathrm{r}_{13},\mathrm{r}_{21}]-[\mathrm{r}_{13},\mathrm{r}_{23}],\]where
\begin{eqnarray*}
\mathrm{r}_{13}.\mathrm{r}_{12}&=&\sum_{i,j}a_i.a_j\otimes b_j\otimes b_i,\quad
\mathrm{r}_{23}.\mathrm{r}_{21}=\sum_{i,j}b_j\otimes a_i.a_j\otimes b_i,\\
\;[\mathrm{r}_{23},\mathrm{r}_{12}]&=&\sum_{i,j}a_j\otimes [a_i,b_j]\otimes b_i,\quad 
\;[\mathrm{r}_{13},\mathrm{r}_{21}]=\sum_{i,j}[a_i,b_j]\otimes a_j\otimes b_i,\quad
\;[\mathrm{r}_{13},\mathrm{r}_{23}]=\sum_{i,j}a_i\otimes a_j\otimes[b_i,b_j].
\end{eqnarray*}

\begin{pr}\label{primportant} For any $\al,\be,\ga\in U^*$, we have
\[ [[\mathrm{r},\mathrm{r}]](\al,\be,\ga)=\prec\ga,\De(\mathrm{r})(\al,\be)\succ. \]
\end{pr}
{\bf Proof.} Recall that according to \eqref{eq110bis}, for any $X\in U$,
\[ \prec[\al,\be],X\succ=\prec\mathrm{L}_{\mathrm{r}_\#(\be)}^t\al
-\mathrm{L}_{\mathrm{r}_\#(\al)}^t\be,X\succ+2\mathrm{L}_X\mathfrak{a}(\al,\be),
 \]$\mathfrak{a}=\frac12\sum_{i}(a_i\otimes b_i-b_i\otimes a_i)$ is the skew-symmetric part of $\mathrm{r}$. 
 We have
\[ \mathrm{r}_\#(\al)=\sum_{i}\prec \al,a_i\succ b_i\esp 
\mathrm{r}_\#(\be)=\sum_{i}\prec \be,a_i\succ b_i. \]So
\begin{eqnarray*}
-\prec\ga,[
\mathrm{r}_\#(\al),\mathrm{r}_\#(\be)]\succ&=&-\sum_{i,j}\prec \al,a_i\succ
\prec \be,a_j\succ\prec \ga,[b_i,b_j]\succ=-[\mathrm{r}_{13},\mathrm{r}_{23}](\al,\be,\ga).\end{eqnarray*}
Now\begin{eqnarray*}
\mathrm{r}(\mathrm{L}_{\mathrm{r}_\#(\be)}^t\al,\ga)&=&\sum_{j}\prec \be,a_j\succ
\mathrm{r}(\mathrm{L}_{b_j}^t\al,\ga)
=\sum_{i,j}\prec \be,a_j\succ
\prec\mathrm{L}_{b_j}^t\al,a_i\succ\prec\ga,b_i\succ\\
&=&\sum_{i,j}(b_j.a_i)\otimes a_j\otimes b_i(\al,\be,\ga).
\end{eqnarray*}
In the same way, we get
\begin{eqnarray*}
-\mathrm{r}(\mathrm{L}_{\mathrm{r}_\#(\al)}^t\be,\ga)&=&-\sum_{i,j}a_j\otimes (b_j.a_i)\otimes b_i(\al,\be,\ga).
\end{eqnarray*}
On the other hand, 
\begin{eqnarray*}
2\mathrm{r}(\mathrm{L}\mathfrak{a}(\al,\be),\ga)&=&2\sum_{i}\prec
\mathrm{L}\mathfrak{a}(\al,\be),a_i\succ\prec \ga,b_i\succ\\
&=&2\sum_{i}\mathfrak{a}(\mathrm{L}_{a_i}^t\al,\be)\prec \ga,b_i\succ+
2\sum_{i}\mathfrak{a}(\al,\mathrm{L}_{a_i}^t\be)\prec \ga,b_i\succ,\\
2\sum_{i}\mathfrak{a}(\mathrm{L}_{a_i}^t\al,\be)\prec \ga,b_i\succ&=&
\sum_{i,j}\left(\prec\mathrm{L}_{a_i}^t\al,a_j\succ\prec\be,b_j\succ
-\prec\mathrm{L}_{a_i}^t\al,b_j\succ\prec\be,a_j\succ
\right)
\prec \ga,b_i\succ\\
&=&\sum_{i,j}\left((a_i.a_j)\otimes b_j\otimes b_i
-(a_i.b_j)\otimes a_j\otimes b_i \right)(\al,\be,\ga),\\
&=&\mathrm{r}_{13}.\mathrm{r}_{12}(\al,\be,\ga)-
\sum_{i,j}(a_i.b_j)\otimes a_j\otimes b_i (\al,\be,\ga),\\
2\sum_{i}\mathfrak{a}(\al,\mathrm{L}_{a_i}^t\be)\prec \ga,b_i\succ&=&
\sum_{i,j}\left(a_j\otimes (a_i.b_j)\otimes b_i
-b_j\otimes (a_i.a_j)\otimes b_i \right)(\al,\be,\ga)\\
&=&\sum_{i,j}a_j\otimes (a_i.b_j)\otimes b_i(\al,\be,\ga)-
\mathrm{r}_{23}.\mathrm{r}_{21}(\al,\be,\ga).
\end{eqnarray*}
By combining all what above we get the desired formula.\hfill$\square$

\begin{remark}\label{rem0}
\begin{enumerate}\item This proposition shows that statement of Theorem \ref{theo11} is the same as the one of Theorem 5.4 in \cite{bai}. However, our proof is more easier because the expression of $\De(\mathrm{r})$ is more simple to handle than  the one of $[[\mathrm{r},\mathrm{r}]]$. The practical nature of $\De(\mathrm{r})$  will be crucial later, in particular, in Sections \ref{section5}-\ref{section5bis}.
\item Let $(U,.)$ be a left symmetric algebra and $\mathrm{r}\in U\otimes U$. According to Proposition \ref{primportant}, $[[r,r]]=0$ iff $\mathrm{r}_\#$ is a Lie algebra endomorphism. This generalizes Theorem 6.6 in \cite{bai}, stated in the case when $\mathrm{r}$ is symmetric.\\
Suppose now that   $\mathrm{r}$ is symmetric and $\mathrm{r}_\#$ is an isomorphism. One can see easily by using \eqref{eq110bis} that, for any $X,Y,Z\in U$, 
\[  \prec\mathrm{r}_\#^{-1}(Z),\De(\mathrm{r})(\mathrm{r}_\#^{-1}(X),\mathrm{r}_\#^{-1}(Y))=
B(X,Y.Z)-B(Y,X.Z)-B(Z,[X,Y]),  \]where $B\in U^*\otimes U^*$ is given by
$B(X,Y)=\prec\mathrm{r}_\#^{-1}(X),Y\succ$. So $[[r,r]]=0$ iff $B$ is 2-cocycle of $(U,.)$. This prove in a different way Theorem 6.3 in \cite{bai}.
\end{enumerate}

\end{remark}

Let introduce now a key notion in our work, namely,  the notion of quasi $S$-matrix as a generalization of the one of $S$-matrix appeared first in \cite{bai}.\\
Let $U$ be a left symmetric  algebra. A \emph{quasi $S$-matrix} of $U$ is a $\mathrm{r}\in U\otimes
U$ such that its skew-symmetric part is $\mathrm{L}_U$-invariant and $[[\mathrm{r},\mathrm{r}]]$ is $Q$-invariant.
 Recall that a \emph{$S$-matrix} of $U$
is a $\mathrm{r}\in U\otimes
U$ which is symmetric and satisfying
\begin{equation}\label{smatrix} [[r,r]]=0.\end{equation}
In what follows, we  focus our attention on the Lie algebra structure on $\Phi(U)$ associated to a quasi $S$-matrix. We show that Lie algebra can be described in a precise and simple way. Indeed, let $\mathrm{r}$  be a quasi $S$-matrix.
 Then, according to Theorem  \ref{theo11},
 the product on $U^*$ given by \eqref{eq110} is left symmetric and
 $(\Phi(U),[\;,\;]^r,\prs_0,K_0)$ is a para-K\"ahler Lie algebra, where
 \[ [X+\al,Y+\be]^r=[X,Y]-\mathrm{L}_{X}^t\be-\mathrm{L}_{\al}^tY+
 \mathrm{L}_{Y}^t\al+\mathrm{L}_{\be}^tX+[\al,\be]. \]  We have shown in Example \ref{exemple1} that $\Phi(U)$ carries a left symmetric product $\triangleright$ and its associated Lie bracket $[\;,\;]^\triangleright$ induces on $\Phi(U)$ a para-K\"ahler Lie algebra structure.
 We  define  a new bracket on $\Phi(U)$ by putting
 \begin{equation}\label{eqnb} [X+\al,Y+\be]^{\triangleright,r}= [X+\al,Y+\be]^\triangleright+\De(\mathrm{r})(\al,\be).\end{equation}
  The following proposition has been inspired to us by a result appeared in \cite{diatta} in the context of Lie bialgebras and $R$-matrices (see Proposition 4.2.1.1 of \cite{diatta}).

\begin{pr}\label{pralg} $(\Phi(U),[\;,\;]^{\triangleright,r})$ is a Lie algebra and
 the linear map $\xi:(\Phi(U),[\;,\;]^{\triangleright,r})\too (\Phi(U),[\;,\;]^r)$, $X+\al\mapsto X-\mathrm{r}_\#(\al)+\al$ is an isomorphism of Lie algebras.
\end{pr}

{\bf Proof.}  Clearly $\xi$ is bijective. Let us show that $\xi$ preserves the Lie brackets. 
It is clear that, for any $X,Y\in U$, $\xi\left([X,Y]^{\triangleright,r} \right)=[\xi(X),\xi(Y)]^r$. Now, for any $X\in U$, $\al\in U^*$,
\begin{eqnarray*}
\xi\left([X,\al]^{\triangleright,r} \right)&=&\xi(-\mathrm{L}_X^t\al)\\
&=&\mathrm{r}_\#(\mathrm{L}_X^t\al)-\mathrm{L}_X^t\al\\
&\stackrel{\eqref{eq110bis}}=&\mathrm{L}_\al^tX-[X,\mathrm{r}_\#(\al)]-\mathrm{L}_X^t\al\\
&=&[X,-\mathrm{r}_\#(\al)+\al]^r\\
&=&[\xi(X),\xi(\al)]^r.
\end{eqnarray*}On the other hand, for any $\al,\be\in U^*$,
\begin{eqnarray*}
\xi\left([\al,\be]^{\triangleright,r} \right)&=&\xi(\De(\mathrm{r})(\al,\be))\\
&=&\De(\mathrm{r})(\al,\be),\\
\;[\xi(\al),\xi(\be)]^r&=&[-\mathrm{r}_\#(\al)+\al,-\mathrm{r}_\#(\be)+\be]^r\\
&=&[\mathrm{r}_\#(\al),\mathrm{r}_\#(\be)]+[\al,\be]+\mathrm{L}_{\mathrm{r}_\#(\al)}^t\be-\mathrm{L}_{\mathrm{r}_\#(\be)}^t\al+\mathrm{L}_{\al}^t\mathrm{r}_\#(\be)-
\mathrm{L}_{\be}^t\mathrm{r}_\#(\al)\\
&\stackrel{\eqref{eq110bis}}=&[\mathrm{r}_\#(\al),\mathrm{r}_\#(\be)]+
\mathrm{r}_\#(\mathrm{L}_{\mathrm{r}_\#(\be)}^t\al)-\mathrm{r}_\#(\mathrm{L}_{\mathrm{r}_\#(\al)}^t\be)+[\mathrm{r}_\#(\be),\mathrm{r}_\#(\al)]-[\mathrm{r}_\#(\al),\mathrm{r}_\#(\be)]\\
&=&\mathrm{r}_\#([\al,\be])-[\mathrm{r}_\#(\al),\mathrm{r}_\#(\be)]\\
&=&\De(\mathrm{r})(\al,\be).
\end{eqnarray*}
\hfill$\square$

We can now transport the para-K\"ahler structure associated to $\mathrm{r}$ from $(\Phi(U),[\;,\;]^r,\prs_0,K_0)$ to $\Phi(U)$ via $\xi$ and we get the following proposition.
\begin{pr}\label{prm} Let $(U,.)$ be a left symmetric algebra and $\mathrm{r}=\mathfrak{a}+\mathfrak{s}\in U\otimes U$ a quasi $S$-matrix. Then $(\Phi(U),[\;,\;]^{\triangleright,r},\prs_r,K_r)$ is a para-K\"ahler Lie algebra, where
\[ \langle X+\al,Y+\be\rangle_r=\langle X+\al,Y+\be\rangle_0-2\mathfrak{s}(\al,\be)\esp K_r(X+\al)=K_0(X+\al)-2\mathrm{r}_\#(\al). \]

\end{pr}

\begin{remark} \label{remdiatta}\begin{enumerate}\item[$(a)$]Actually, by using a similar method, we can generalize the result of Diatta \cite{diatta}. Let $(\G,[\;,\;])$ be a Lie algebra and $\mathrm{r}\in\G\wedge\G$. Define on $\G^*$ and $\Phi(\G)$, respectively, two brackets $[\;,\;]^*$ and $[\;,\;]^\mathrm{r}$ by
\[ [\al,\be]^*=\ad^*_{\mathrm{r}_\#(\al)}\be 
-\ad^*_{\mathrm{r}_\#(\be)}\al
\esp[X+\al,Y+\be]^\mathrm{r}=[X,Y]+[\al,\be]^*-\ad_X^t\be-\ad_\al^tY
+\ad_Y^t\al+\ad_\be^tX,\]and $[\mathrm{r},\mathrm{r}]\in\G\otimes\G\otimes\G\simeq\mathrm{End}(\G^*\otimes\G^*,\G)$ by
\[ [\mathrm{r},\mathrm{r}](\al,\be)=\mathrm{r}_\#([\al,\be]^*)-[\mathrm{r}_\#(\al),\mathrm{r}_\#(\be)]
. \]
It is well-known that $[\;,\;]^*$ is a Lie bracket iff 
$[\mathrm{r},\mathrm{r}]$ is $\ad$-invariant. In this case, $[\;,\;]^\mathrm{r}$ is a Lie bracket.
 Define a new bracket on $\Phi(\G)$ by putting
\[ [X+\al,Y+\be]^{\diamond,r}=[X,Y]+\ad_X^*\be-\ad_Y^*\al+[\mathrm{r},\mathrm{r}](\al,\be). \]
 By using the same argument in the proof of Proposition \ref{pralg}, one can see that $(\Phi(\G),[\;,\;]^{\diamond,r})$ is a Lie algebra and
 the linear map $\xi:(\Phi(\G),[\;,\;]^{\diamond,r})\too (\Phi(\G),[\;,\;]^r)$, $X+\al\mapsto X-\mathrm{r}_\#(\al)+\al$ is an isomorphism of Lie algebras. When  $[\mathrm{r},\mathrm{r}]=0$, we recover the result of Diatta.
 \item[$(b)$] Let $(U,.)$ be a left symmetric algebra and $\mathrm{r}=\mathfrak{a}+\mathfrak{s}\in U\otimes U$ a quasi $S$-matrix. The Lie algebra $(\Phi(U),[\;,\;]^{\triangleright,r})$ is a $\Z_2$-graded Lie algebra and hence $L: U^*\times U^*\times U^*\too U^*$ given by
 \[ L(\al,\be,\ga)=\mathrm{L}_{\De(\mathrm{r})(\al,\be)}^*\ga \]
 is a Lie triple system (see for instance \cite{jacobson, lister, smirnov} for the definition and the properties of Lie triple systems).\end{enumerate}
\end{remark}

\section{Some classes of para-K\"ahler  Lie
algebras}\label{section4}

In this section, we develop some methods to build para-K\"ahler Lie algebras
 based on  the following
proposition where we adopt the notations of the last section, in particular,  Proposition \ref{prm}.
\begin{pr}\label{pr31}Let $(U,.)$ be a left symmetric algebra and $\mathrm{r}=
\mathfrak{a}+\mathfrak{s}\in U\otimes
U$ which is 
$\mathrm{L}_U$-invariant. Then
$\mathrm{L}(\mathfrak{a})=0$,  $[[\mathrm{r},\mathrm{r}]]=\De(\mathrm{r})=0$ and $(\Phi(U),[\;,\;,]^{\triangleright},\prs_r,K_r)$  a para-K\"ahler Lie algebra. Moreover, the Levi-Civita product of $(\Phi(U),[\;,\;,]^{\triangleright},\prs_r)$ is $\triangleright$ given by \eqref{exem1eq1}.

\end{pr}

{\bf Proof.} Since $\mathrm{r}$ is $\mathrm{L}_U$-invariant then so $\mathfrak{a}$  and
hence $\mathrm{L}(\mathfrak{a})=0$. The vanishing of $\De(\mathrm{r})$ is immediate. So we can apply Proposition \ref{prm}. To conclude, one can check easily that $\triangleright$ is actually the Levi-Civita product of 
$(\Phi(U),[\;,\;,]^{\triangleright},\prs_r)$. 
\hfill$\square$

As consequence of this proposition we get the following large class of para-K\"ahler Lie
algebras.

\begin{pr}\label{pr32} Let $(\G,\prs)$ be a  pseudo-Riemannian flat Lie algebra, ie., the Levi-Civita product $"."$ is left symmetric. Denote by
$\flat:\G\too \G^*$  the isomorphism associated to $\prs$. Then  $(\Phi(\G),[\;,\;,]^{\triangleright},\prs_\flat,K_\flat)$ is a para-K\"ahler Lie algebra, where
\[ \langle X+\al,Y+\be\rangle_\flat=\langle X+\al,Y+\be\rangle_0-2\langle\flat^{-1}(\al),\flat^{-1}(\be)\rangle\esp K_\flat(X+\al)=K_0(X+\al)-2\flat^{-1}(\al). \]
Moreover, the Levi-Civita product of $(\Phi(\G),[\;,\;,]^{\triangleright},\prs_\flat)$ is $\triangleright$ given by \eqref{exem1eq1}.

 \end{pr}

{\bf Proof.} The Levi-Civita product defines a left symmetric algebra structure on
$U$ and  $\mathrm{r}\in U\otimes U$ defined by
$\mathrm{r}(\al,\be)=\langle\flat^{-1}(\al),\flat^{-1}(\be)$
is $\mathrm{L}_U$-invariant and one can conclude by using Proposition \ref{pr31}.\hfill $\square$\\

Let us give some methods to build pseudo-Riemannian flat Lie algebras.\\Let $({
U},[.,.],\om)$ be a  symplectic Lie algebra and $\mathrm{B}$ a nondegenerate bi-invariant
bilinear
symmetric form on $U$.  The isomorphism $D$ defined by
$$\om(X,Y)=
\mathrm{B}(D(X),Y)$$ is an invertible derivation and hence $U$ is nilpotent (see
\cite{Jacob}).
 The nondegenerate symmetric
bilinear form $\prs$ given by
$$\langle X,Y\rangle=\mathrm{B}(D(X),D(Y))$$satisfies
$$\langle X.Y,Z\rangle+\langle Y,X.Z\rangle=0,$$ where the dot designs the left symmetric
product associated to $\om$ given by \eqref{eq15}.  Thus  $(U,\prs)$ is a flat
pseudo-Riemannian Lie
algebra (see \cite{BBM}). So any symplectic quadratic Lie algebra $(U,\mathrm{B},\om)$
gives rise to
a flat pseudo-Riemannian Lie
algebra $(U,\prs)$. \\ More generally, let $(\G,[\;,\;],B)$ be a quadratic Lie algebra
and $\mathfrak{r}\in\G\wedge\G$ a solution of the classical Yang-Baxter equation. The product on $\G^*$ given by $\al.\be=\ad_{\mathfrak{r}_\#(\al)}^*\be$ is left
symmetric and hence induces a Lie bracket $[\;,\;]_{\mathfrak{r}}$ on $\G^*$. In fact,
this product is the Levi-Civita product of $B^*$ (the induced bilinear nondegenerate
symmetric form on $\G^*$). Thus $(\G^*,[\;,\;]_{\mathfrak{r}},B^*)$ is a flat
pseudo-Riemannian Lie algebra (see \cite{bm}).\\
Let us give now a method to built symplectic quadratic Lie algebras.\\
Let $n \in {\Bbb N}^*$ and 
${\mathcal  A}$ a vector space with a basis $\{e_1,\dots,e_n\}.$    We
consider on ${\mathcal  A}$ the
product
defined by
$$e_ie_j= e_je_i= e_{i+j}\,\, \, \mbox{if} \,\,  i+j\leqq n,\,\,\,  e_ie_j= e_je_i=
0\,\,\,  
\mbox{if}\,\,  i+j>n.$$
The vector space ${\mathcal  A}$  endowed with this product is a commutative and
associative algebra.
\\Let $({\mathcal L},[.,])$ be an arbitrary Lie algebra. Then the following product
$$[X\otimes a,Y\otimes b]_{\mathcal T}:= [X,Y]\otimes ab,$$
defines a structure of Lie algebra on the vector space ${\mathcal T}:= {\mathcal L}\otimes
{\mathcal 
A}$. Moreover
the 
endomorphism $\delta$ of ${\mathcal T}$ defined by
$$\delta(X\otimes e_i):= i X\otimes e_i$$ for any $
X\in {\mathcal L}$ and any $i\in\{1,\dots,n\},$ is an invertible derivation of
${\mathcal T}$.\\
Now, on the vector space ${ U}:= {\mathcal T}\oplus {\mathcal T}^*$ we will
define a
structure of symplectic quadratic algebra in the following way. For any $s,t\in{\mathcal
T}$ and
any $f,h\in{\mathcal T}^*$, put
\begin{eqnarray*}
 [t+f,s+h]_{ U}&=&  [t,s]_{\mathcal T}-h\circ{\ad}_{\mathcal
T}(t)+f\circ{\ad}_{\mathcal
T}(s),\\
\mathrm{B}(t+f,s+h)&=& f(s)+h(t),\\
D(t+f)&=& \delta(t)-f\circ \delta,\\
\om(t+f,s+h)&=&\mathrm{B}(D(t+f),s+h).
\end{eqnarray*}
One can check easily that $(U,\mathrm{B},\om)$ is a symplectic quadratic Lie algebra.

\begin{pr} \label{pr34} Let $(\G,[\;,\;])$ be a Lie algebra,  $\mathrm{b}\in \wedge^2
\G$ is a solution
of the
classical Yang-Baxter equation on $(\G,[\;,\;])$, i.e. $[\mathrm{b},\mathrm{b}]=0$, and  $\mathrm{r}=\mathfrak{s}+\mathfrak{a}\in
\G^*\otimes \G^*$ such that $
 \ad_{\mathrm{b}_\#(\al)}^*\mathrm{r}=0,
$ for any $\al\in \G^*$.
Then    $(\Phi(\G),[\;,\;]^\mathrm{b},\prs^r,K^r)$ is a para-K\"ahler Lie algebra, where
\[ [X+\al,Y+\be]^\mathrm{b}=\ad^*_{\mathrm{b}_\#(\al)}\be-\ad^*_{\mathrm{b}_\#(\be)}\al+[\mathrm{b}_\#(\al),Y]+[X,\mathrm{b}_\#(\be)], \]
\[ \langle X+\al,Y+\be\rangle^r=\langle X+\al,Y+\be\rangle_0-2\mathfrak{s}(X,Y)\esp K^r(X+\al)=-K_0(X+\al)-2\mathrm{r}_\#(X). \]
Moreover, the Levi-Civita product associated to $(\Phi(\G),[\;,\;]^\mathrm{b},\prs^r)$ is left symmetric and it is given by
\[ (X+\al)\triangleright_{\mathrm{b}} (Y+\be)=
\ad^*_{\mathrm{b}_\#(\al)}\be+[\mathrm{b}_\#(\al),Y].\]

\end{pr}
 
 {\bf Proof.} It is well-known that the product on $\G^*$ given by
 $\al.\be=\ad_{\mathrm{b}_\#(\al)}^*\be$ is left symmetric and the condition
$ \ad_{\mathrm{b}_\#(\al)}^*\mathrm{r}=0$ is equivalent to $\mathrm{r}$ is invariant with
respect to this product on
$\G^*$. So  $(\G^*,.,\mathrm{r})$ satisfies the hypothesis of Proposition \ref{pr31}
and the proposition follows.\hfill$\square$

 There is an interesting case of this situation.
 \begin{co}\label{co31} Let $(\G,[\;,\;])$ be a  Lie
algebra, 
$\mathrm{b}\in \wedge^2 \G$ is a solution of the
classical Yang-Baxter equation on $(\G,[\;,\;])$ and $\mathrm{k}\in \G^*\otimes \G^*$
is the Killing form. Then   $(\Phi(\G),[\;,\;]^\mathrm{b},\prs^\mathrm{k},K^\mathrm{k})$ is a para-K\"ahler Lie algebra.\end{co}

\section{Hyper-para-K\"ahler Lie algebras}\label{section3}
 Hyper-para-K\"ahler Lie algebras known also as hyper-symplectic Lie algebras constitute a subclass of the class of para-K\"ahler Lie
algebras. We
will use our study in the last sections to give a new characterization of these Lie algebras. This characterization leads to a notion of compatibility between two left symmetric algebra structures on a given vector space. Since a hyper-para-K\"ahler Lie algebra has a complex product structure we get also a characterization of such structures.
\\

A \emph{hyper-para-K\"ahler Lie algebra} is   a para-K\"ahler Lie algebra
$(\G,\prs,K)$
endowed with an endomorphism $J$ such that, $J^2=-\mathrm{Id}_{\G}$, $JK=-KJ$, $J$ is
skew-symmetric with respect to $\prs$ and $J$ is invariant with respect to the Levi-Civita
product. According to Theorem \ref{theoextendible}, a para-K\"ahler Lie algebra can be identified to the phase space of Lie-extendible left symmetric algebras so it is natural to see how hyper-para-K\"ahler Lie algebras can be described in this sitting.
\begin{pr}
\label{pr31b} Let $(U,.)$ and $(U^*,.)$ be a couple of Lie-extendible left symmetric algebras,  $(\Phi(U),\prs_0,K_0)$ the associated para-K\"ahler Lie algebra and $J:\Phi(U)\too \Phi(U)$ an endomorphism. Then $(\Phi(U),\prs_0,K_0,J)$ is a hyper-para-K\"ahler Lie algebra if and only if there exists a bilinear nondegenerate $\om\in\wedge^2U^*$ such that:
\begin{enumerate}\item[$(i)$] for any $X\in U,\al\in U^*$, $JX=\flat(X)$ and $J\al=-\flat^{-1}(\al)$ where $\flat:U\too U^*$ is the  isomorphism  given by $\flat(X)=\om(X,.)$,
\item[$(ii)$] $(U,.,\om)$ and $(U,\circ,\om)$ are symplectic left symmetric algebras where $\circ$ is given by
$$X\circ Y=\flat^{-1}(\flat(X).\flat(Y)).$$

\end{enumerate}
\end{pr}

{\bf Proof.} From
the relation $JK_0=-K_0J$, we deduce that for any $X\in U$, $JX\in U^*$ and hence $J$ defines
an
isomorphism $\flat:U\too U^*$. Moreover, from $J^2=-\mathrm{Id}_{\G}$ we deduce
that
$J\al=-\flat^{-1}\al$ for any $\al\in U^*$.  The skew-symmetry of $J$ implies that
 $\om\in\wedge^2 U^*$ given by $$\om(X,Y)=\prec\flat(X),Y\succ$$is skew-symmetric and, actually is nondegenerate. 
Now $J$ is invariant if and only if, for any $X,Y,\in U$ and any $\al,\be\in U^*$, 
\[ -\mathrm{L}_X^t(JY)=J(X.Y),\; X.J(\al)=-J(\mathrm{L}_X^t\al),\; \al.J(X)=-J(\mathrm{L}_\al^tX) 
\esp -\mathrm{L}_\al^t(J\be)=J(\al.\be). \]
This is equivalent to
\[ \mathrm{L}_X^t\circ\flat+\flat\circ \mathrm{L}_X=
\flat^{-1}\circ \mathrm{L}_X^t+\mathrm{L}_X\circ\flat^{-1}=0 \esp
 \mathrm{L}_\al^t\circ\flat^{-1}+\flat^{-1}\circ \mathrm{L}_\al=
 \flat\circ \mathrm{L}_\al^t+\mathrm{L}_\al\circ\flat=0,\]for any $X\in U,\al\in U^*$. Now it is obvious that these relations are equivalent to
\begin{equation}\label{eqcom} \mathrm{L}_X^t\circ\flat+\flat\circ \mathrm{L}_X=
0 \esp
 \flat\circ \mathrm{L}_\al^t+\mathrm{L}_\al\circ\flat=0,\end{equation}for any $X\in U,\al\in U^*$. One can see easily that this is equivalent to
$$\om(X. Y,Z)+\om( Y,X. Z)=0\esp\om(X\circ Y,Z)+\om( Y,X\circ Z)=0,
$$for any $X,Y,Z\in U$. 
Thus $(U,.,\om)$ and $(U,\circ,\om)$ are symplectic left symmetric algebras. The converse is obviously true.\hfill$\square$\\

Now, let $U$ be a vector space,  $\om\in\wedge^2U^*$ nondegenerate and $\bullet,\circ$  two  products on $U$ such that $(U,\bullet,\om)$ and $(U,\circ,\om)$ are symplectic left symmetric algebras.  Define $J_0:\Phi(U)\too\Phi(U)$ by $J_0X=\flat(X)$ and $J_0\al=-\flat^{-1}(\al)$ and denote by $\flat(\circ)$ the product on $U^*$ image by $\flat$ of $\circ$. Let the dot denote the product on $\Phi(U)$ extending $(U,\bullet)$ and $(U^*,\flat(\circ))$ by \eqref{eq16}. It is easy to see, by using \eqref{eqcom}, that for any 
$X,Y\in U$, $\al,\be\in
U^*$, 
\begin{equation}\label{eq22}
(X+\al).(Y+\be)=X\bullet Y+\flat^{-1}(\al)\circ
Y-\left({\mathrm{L}}_{\flat^{-1}(\al)}^\circ\right)^t\be
-\left(\mathrm{L}_X^\bullet\right)^t\be.
\end{equation}
We have shown in Proposition \ref{pr13} that this product is Lie-admissible iff \eqref{eq19} hold. Thank to $\om$ we can identify $\Phi(U)$ to $T(U)$. Indeed, 
define $\xi:T(U)\too\Phi(U)$ by
\[ \xi(X,0)=X\esp\xi(0,X)=\flat(X). \]
We have $\Om_1=\xi^*\Om_0,$ $\prs_1=\xi^*\prs_0$, $K_1=\xi^{-1}\circ K_0\circ\xi$ and $J_1=\xi^{-1}\circ J_0\circ\xi$.
It is easy to check that
\begin{equation}\label{eq22bis} (X,Y).(Z,T):=\xi^{-1}(\xi(X,Y).\xi(Z,T))=(X\bullet Z,X\bullet T)+(Y\circ Z,Y\circ T).
\end{equation}
Now on can see easily by using \eqref{eqcom} that, for any $X\in U$ and any $\al\in U^*$, 
\begin{eqnarray*}
 \rho(X,\al)&=&-
\mathrm{K}^{\bullet,\circ}(X,\flat^{-1}
(\al))\esp
\rho^*(\al,X)=\flat\circ \mathrm{K}^{\bullet,\circ}(X,\flat^{-1}(\al))\circ\flat^{-1}, 
\end{eqnarray*}where
$$
 \mathrm{K}^{\bullet,\circ}(X,Y)=[\mathrm{L}_{X}^\bullet,{\mathrm{L}}_Y^\circ]
-\left(\mathrm {L }_{
{X}\bullet Y}^\circ-\mathrm{L}_{Y\circ X}^\bullet\right).
$$(To distinguish between $\bullet$ and $\circ$, we denote by $\mathrm{L}^\bullet_X$ the left multiplication by $X$ associated to $\bullet$ and so on). So by using Proposition \ref{pr13}, we get the following proposition which, actually, does not involves $\om$.

 \begin{pr}\label{pr21}Let $U$ be a vector space and $\bullet,\circ$  two left symmetric  products on $U$. The following assertions are equivalent:
\begin{enumerate}\item
 The product given by \eqref{eq22bis} is Lie-admissible.  
\item For any $X,Y,Z\in U$, $
\mathrm{K}^{\bullet,\circ}(X,Y)Z=\mathrm{K}^{\bullet,\circ}(Z,Y)X$ and
 $
\mathrm{K}^{\bullet,\circ}(X,Y)Z=\mathrm{K}^{\bullet,\circ}(X,Z)Y.$
  
\end{enumerate}

Moreover,, the  product given by \eqref{eq22bis} is left symmetric if and only if 
$\mathrm{K}^{\bullet,\circ}$ vanishes identically.

\end{pr}One can see easily that the second assertion in this proposition
is equivalent to{\small
\begin{eqnarray}
Y\circ[X,Z]^\bullet-[Y\circ X,Z]^\bullet-[X,Y\circ Z]^\bullet&=&
(Z\bullet Y)\circ X-(X\bullet Y)\circ Z,\label{com1}\\
Y\bullet[X,Z]^\circ-[Y\bullet X,Z]^\circ-[X,Y\bullet Z]^\circ&=&
(Z\circ Y)\bullet X-(X\circ Y)\bullet Z,\label{com2}
\end{eqnarray}}for any $X,Y,Z\in U$.

Let us state an important formula.  Let $\bullet$ and $\circ$ be two algebra structures on a vector space $U$.   A straightforward computation gives the following formula:
\begin{equation}\label{eq21}
 {\mathrm{K}^{\bullet+\circ}}(X,Y)=\mathrm{K}^{\bullet}(X,Y)+\mathrm{K}^{\circ}(X,Y)+\mathrm{K}^{\bullet,\circ}(X,Y)-
 \mathrm{K}^{\bullet,\circ}(Y , X),
\end{equation}where $\mathrm{K}^{x}$ is the curvature of $x$.

\begin{Def}\label{Def1}

Two   left symmetric algebras structures $\bullet$ and $\circ$ on $U$ will be called compatible if they satisfy \eqref{com1}-\eqref{com2}
or equivalently
$\mathrm{K}^{\bullet,\circ}$ satisfies the second assertion in Proposition \ref{pr21}.\end{Def}

  The following proposition is an immediate consequence of \eqref{eq21}.
\begin{pr}\label{prcompatible} Let $\bullet,\circ$ be two compatible left symmetric algebra structures on $U$. Then for any $a,b\in\R$,\\ $(U,a\;\bullet+b\;\circ)$ is a left symmetric algebra.

\end{pr}

\begin{remark} Let $\bullet,\circ$ be two compatible left symmetric algebra structures on $U$. As consequence of Proposition \ref{prcompatible}, the bracket $a[\;,\;]^\bullet+b[\;,\;]^\circ$ is Lie bracket and hence the two dual Poisson structures on $U^*$ associated to $[\;,\;]^\bullet$ and $[\;,\;]^\circ$ are compatible (see for instance \cite{vaisman} for the definition of compatible Poisson structures).

\end{remark}
Finally, we get a characterization of  hyper-para-K\"ahler Lie algebras. Actually, our study can be generalized easily to give a characterization of complex product structures.  The characterization given in the following theorem completes the study of complex product structures achieved in \cite{salomon}.

\begin{theo}\label{theohyperkahler}\begin{enumerate}\item Let $\bullet,\circ$ be two compatible left symmetric algebra structures on $U$. Then $(T(U),K_1,J_1)$ endowed with the Lie algebra structure associated to the product given by \eqref{eq22bis} is a complex product Lie algebra. Moreover, all complex product Lie algebras are obtained in this way.
\item Let $\bullet,\circ$ be two compatible left symmetric algebra structures on $U$ and $\om\in\wedge^2U^*$ such that $(U,\om,\bullet)$ and
$(U,\om,\circ)$ are symplectic left symmetric algebras. Then $(T(U),\prs_1,K_1,J_1)$ endowed with the Lie algebra structure associated to the product given by \eqref{eq22bis} is a hyper-para-K\"ahler Lie algebra. Moreover, all hyper-para-K\"ahler Lie algebras are obtained in this way.

\end{enumerate}

\end{theo}

{\bf Proof.}\begin{enumerate}\item To show that $(T(U),K_1,J_1)$ is a complex product Lie algebra it suffices to show that the Nijenhuis torsion of $K_1$ and $J_1$ vanishes which is easy to check. Conversely, let $(\G,K,J)$ be a complex product Lie algebra. We have $\G=\G^1\oplus\G^{-1}$ where $\G^i=\ker(K+i\mathrm{Id}_\G)$ and $J$ defines an isomorphism $\phi:\G^1\too\G^{-1}$. 
We consider  the product
$"."$ on $\G$ given by
\[ (u^1+u^{-1}).(v^1+v^{-1})=u^1\circ v^1+\phi(u^1\circ\phi^{-1}(v^{-1}))+
 \phi^{-1}(u^{-1}\star\phi (v^{1}))+u^{-1}\star v^{-1},\]where $\circ$ and $\star$ are the products on $\G^1$ and $\G^{-1}$, respectively, given by
 $$u^1\circ v^1=-\pi_1J[u^1,Jv^1]\esp u^{-1}\star v^{-1}=-\pi_{-1}J[u^{-1},Jv^{-1}],$$where $\pi_i$ is the projection on $\G^i$.
 It was shown in \cite{salomon} that $\circ,\star$ are left symmetric and $"."$ is Lie-admissible. Put $U=\G^1$, $\bullet=\phi^{-1}(\star)$ and define $\xi:T(U)\too \G$ by $\xi(X,0)=X$ and $\xi(0,X)=\phi(X)$. We get the desired isomorphism.
 \item This is a consequence of the study above.\hfill$\square$

\end{enumerate}

\begin{exem}
Let $(U,.)$ be a left symmetric algebra. Then $"."$ is compatible with itself so $(T(U),K_1,J_1)$ endowed with the Lie algebra bracket associated to the left symmetric product
\[ (X,Y).(Z,T)=(X. Z+Y. Z,Y. T+X. T) \]is a complex product Lie algebra. Moreover, if $U$ carries $\om$ such that $(U,.,\om)$
is a symplectic left symmetric algebra then $(T(U),\prs_1,K_1,J_1)$ is a hyper-para-K\"ahler Lie algebra.
\end{exem}

The following proposition is immediate.
\begin{pr}\label{prabelian} Let $\bullet,\circ$ be two compatible left symmetric algebra structures on $U$ and $(T(U),K_1,J_1)$ the associated complex product structure. Then the following are equivalent:
\begin{enumerate}\item[$(i)$] $K_1$ is abelian.
\item[$(ii)$] $J_1$ is abelian.
\item[$(iii)$] Both $\bullet$ and $\circ$ are commutative and hence associative.

\end{enumerate}

\end{pr}

According to  \eqref{com1}-\eqref{com2}, two associative and commutative algebra structures $\bullet$ and $\circ$ on $U$ are compatible if, for any $X,Y,Z\in U$,
\begin{eqnarray}
(Z\bullet Y)\circ X-(X\bullet Y)\circ Z&=&
(Z\circ Y)\bullet X-(X\circ Y)\bullet Z=0.\label{com3}
\end{eqnarray}
In this case, $(T(U),K_1,J_1)$ endowed with the bracket associated to the product given by \eqref{eq22bis} is an abelian complex product structure. There is similar result in \cite{salomon} with a different product.

\section{Quasi $S$-matrices on symplectic Lie algebras}\label{section5}
We have shown in Section \ref{section2} that finding the set of quasi $S$-matrices on a given left symmetric algebra gives a large class of para-K\"ahler Lie algebras. In this section, we investigate the set of quasi $S$-matrices with respect to the left symmetric product associated to a symplectic Lie algebra.

Let $(\G,\om)$ be a symplectic Lie algebra and $\flat:\G\too\G^*$ the isomorphism given by $\flat(X)=\om(X,.)$. The product $\mathbf{a}$ given  by \eqref{eq15} is left symmetric. We associate to any endomorphism $A:\G\too\G$ the tensor $\mathrm{YB}(A)\in\mathrm{End}(\G\otimes \G,\G)$ given by
\begin{equation}\label{yb} \mathrm{YB}(A)(X,Y)=A[AX,Y]+A[X,AY]-[AX,AY].
\end{equation}
The following proposition gives an useful characterization of quasi $S$-matrices  on $(\G,\mathbf{a})$.
\begin{pr}\label{prqss}Let $\mathrm{r}\in \G\otimes\G$ and $\mathfrak{a}$ its skew-symmetric part. Put $A=\mathrm{r}_\#\circ\flat$ and $T=\mathfrak{a}_\#\circ\flat$. Then the following assertions holds:
\begin{enumerate}\item[$(i)$] The tensor $\mathrm{r}$ is a quasi $S$-matrix of
$(\G,\mathbf{a})$ if and only if $T$ and $\mathrm{YB}(A)$ are $\ad$-invariant.
\item[$(ii)$] If $\mathrm{r}$ is symmetric then it is a  $S$-matrix of
$(\G,\mathbf{a})$ if and only if  $\mathrm{YB}(A)=0$.
\item[$(iii)$] If $\mathrm{r}$ is symmetric and invertible then it is a  $S$-matrix of
$(\G,\mathbf{a})$ if and only if  $A^{-1}$ is a derivation of the Lie algebra $\G$.

 \end{enumerate}
\end{pr}

 {\bf Proof.} From \eqref{eq15}, we get that for any $X\in\G$, \begin{equation}\label{sympl}\flat\circ\ad_X=\mathrm{L}_X^*\circ\flat.
 \end{equation} Moreover, we have, for any $X,Y\in\G$,
 \begin{eqnarray*}\De(\mathrm{r})(\flat X,\flat Y)&=&\mathrm{r}_\#([\flat X,\flat Y])-
 [\mathrm{r}_\#(\flat X),\mathrm{r}_\#(\flat Y)]
 \\&\stackrel{\eqref{eq110bis}}=&
 A\circ\flat^{-1}\left(\mathrm{L}_{AY}^t\flat X \right)
 -A\circ\flat^{-1}\left(\mathrm{L}_{AX}^t\flat Y \right)\\&&+2
 A\circ\flat^{-1}\left(\mathrm{L}(\mathfrak{a})(\flat X,\flat Y) \right)-[AX,AY]\\
 &\stackrel{\eqref{sympl}}=&\mathrm{YB}(A)(X,Y)+2
  A\circ\flat^{-1}\left(\mathrm{L}(\mathfrak{a})(\flat X,\flat Y) \right).
 \end{eqnarray*}
 From this relation and \eqref{sympl} one can see easily that $(i)$ and $(ii)$ holds. Now it is easy that if $A$ is invertible then $\mathrm{YB}(A)=0$ if and only if $A^{-1}$ is a derivation of the Lie algebra $\G$ and $(iii)$ holds.\hfill$\square$

 Let $\G$ be a Lie algebra. The modified Yang-Baxter equation is the equation
 \begin{equation}\label{myb}
 \mathrm{YB}(A)(X,Y)=t[X,Y],\quad\mbox{for all}\; X,Y\in\G,
 \end{equation}where $t\in\R$ is a fixed parameter and the unknown $A$ is an endomorphism of $\G$. When $t=0$ we get the operator form of the classical Yang-Baxter equation. The following proposition is an immediate consequence of Proposition \ref{prqss}.
 \begin{pr}Let $(\G,\om)$ be a symplectic Lie algebra and $A$ a solution of the modified Yang-Baxter equation which is skew-symmetric with respect to $\om$. Then $\mathrm{r}=A\circ\flat^{-1}$ is a quasi $S$-matrix of $(\G,\mathbf{a})$.
 
 \end{pr}

\begin{theo}\label{theosymp}Let $(\G,\om)$ be a symplectic Lie algebra and $A:\G\too\G$. Put $A=A^s+A^a$ where $A^s$ and $A^a$ are, respectively, the symmetric and the skew-symmetric part of $A$ (with respect to $\om$). If  both $\mathrm{YB}(A)$ and $A^s$ are $\ad$-invariant, then the product $\circ$ on $\G$ given by $X\circ Y=X.[(A^s-A^a)Y]-(AX).Y$ is left symmetric and
$(T(\G),\prs_A,K_A)$ endowed with the Lie bracket  given by
\[  [(X,Y),(Z,T)]^A=([X,Z]+\mathrm{YB}(A)(Y,T),[X,T]+[Z,Y] ) \] is a para-K\"ahler Lie algebra, where
\[ \langle(X,Y),(Z,T)\rangle_A=\om(T,X)+\om(Y,Z)
+2\om(A^aY,T)\esp
K_A(X,Y)=(X-2AY,-Y), \]
and the dot is the left symmetric product associated to $(\G,\om)$.

\end{theo}

{\bf Proof.} According to Proposition \ref{prqss}, $\mathrm{r}$ given by 
$\mathrm{r}(\al,\be)=-\om( A\flat^{-1}(\al),\flat^{-1}(\be))$ is a quasi $S$-matrix with respect to the left symmetric product associated to $\om$.  By virtue of  Corollary \ref{co11} and Proposition \ref{prm}, $\mathrm{r}$ defines on
$\G^*$ a left
symmetric Lie algebra structure by \eqref{eq110} and $(\Phi(\G),[\;,\;]^{\triangleright,r}\prs_r,K_r)$  
is a para-K\"ahler Lie algebra.  We consider now the linear isomorphism $\mu:T(\G)\too\Phi(\G)$, $(X,Y)\mapsto(X,\flat(Y))$. Thus $(T(\G),[\;,\;]^\mu,\mu^*\prs_r,
\mu^{-1}\circ K_r\circ \mu)$ is a para-K\"ahler Lie algebra, $[\;,\;]^\mu$ is a pull-back by $\mu$ of $[\;,\;]^{\triangleright,r}$.  One can check easily that this bracket is the Lie bracket given in the statement of the theorem,   $\prs_A=\mu^*\prs_r$ and $K_A=\mu^{-1}\circ K_r\circ \mu$. \\
Let us compute now the pull-back by $\flat$ of the left symmetric product on $\G^*$ given by \eqref{eq110}. We have
\begin{eqnarray*}
\prec\al,X\circ Y\succ&=&-\prec\flat(X).\flat(Y),\flat^{-1}(\al)\succ\\
&\stackrel{\eqref{eq110}}=&-\mathrm{r}(\mathrm{L}_{\flat^{-1}(\al)}^t\flat(X),\flat(Y))
-\mathrm{r}(\flat(X),\ad_{\flat^{-1}(\al)}^t\flat(Y))\\
&\stackrel{\eqref{sympl}}=&\om(A[\flat^{-1}(\al),X],Y)
-\om(AX,\flat^{-1}(\ad_{\flat^{-1}(\al)}^t\flat(Y)))\\
&=&\om([\flat^{-1}(\al),X],A^sY)+\om([X,\flat^{-1}(\al)],A^aY)+\om(Y,[\flat^{-1}(\al),AX])\\
&=&\prec\al,X.[(A^s-A^a)Y]-(AX).Y\succ.
\end{eqnarray*}
\hfill$\square$

\begin{remark}\label{rem3}Actually, this theorem and Remark \ref{remdiatta} $(b)$ suggest as the following more general result. Let $\G$ be a Lie algebra and $A$ an endomorphism of $\G$ such that $\mathrm{YB}(A)$ is $\ad$-invariant. Then, one can check that the bracket $[\;,\;]^A$ on $T(\G)$ is a Lie bracket and hence $L^A:\G\times\G\times\G\too\G$ given by
\[ L^A(X,Y,Z)=[\mathrm{YB}(A)(X,Y),Z] \]is a Lie triple system.

\end{remark}

\begin{exem}  Let $(U,.)$ be a left symmetric algebra with an invertible derivation $D$. We know that on the vector space $\Phi(U):= U\oplus U^*$ we have a left symmetric structure defined by:
$$(X+\alpha)\rhd (Y+\beta):= X.Y - L_X^t\beta,\,\,\,\,\, \forall X,Y \in U, \alpha,\beta \in U^*.$$
 Moreover $(\Phi(U),\rhd,\Gamma_0)$ is a symplectic  left symmetric algebra.
 Now, it is easy to verify that the endomorphism $\Delta$ of $\Phi(U)$ defined by:
$$\Delta(X+\alpha):= D(X) - \alpha\circ D, \forall X \in U, \alpha \in U^*,$$
is an invertible derivation of $(\Phi(U),\rhd)$ which is skew-symmetric with respect to $\Gamma_0.$
 We are going, in the following, to construct a left symmetric algebras with an invertible derivation.
 Let $n\in {\Bbb N}^*$ and $A$ a vector space with a basis $\{e_1,\dots,e_n\}$. We consider on $A$ the product defined by:
$$ e_ie_j= e_je_i:= e_{i+j}\,\,\, \mbox{if}\,\, i+j\leq n, \,\, \, e_ie_j= e_je_i:= 0 \,\,\, \mbox{if}\,\, i+j> n.$$
The vector space $A$ endowed with this product is  a commutative associative algebra.
 Let $(V,\star)$ be a symmetric algebra, then $(U:= V\otimes A,.)$ is a left symmetric algebra where the product "." is defined by:
 $$ v\otimes e_i.w\otimes e_j:= v\star w\otimes e_ie_j,\,\,\, \forall v, w\in V, i,j \in\{1,\dots,n\}.$$
 Moreover, the endomorphism $D$ of $U$ defined by:
 $$D(v\otimes e_i):= i v\otimes e_i,\,\,\, \forall v\in V, i\in \{1,\dots,n\},$$
 is an invertible derivation of $(U,.)$.
 Finally, by using the first construction, we obtain a symplectic left symmetric algebra $(\Phi(U),\rhd,\Gamma_0)$ with an invertible derivation $\Delta$ which is  skew-symmetric with respect to $\Gamma_0.$

\end{exem}

\section{Quasi $S$-matrices on a left symmetric  algebra $U$ with an invariant isomorphism
$\Theta:U\too U^*$}\label{section5bis}

In this section, we investigate the set of quasi $S$-matrices on  a left symmetric  algebra $U$ with an invariant isomorphism
$\Theta:U\too U^*$. The most important classes are symplectic left symmetric algebras and flat pseudo-Riemannian Lie algebras.\\

Let $(U,.)$ be a left symmetric algebra and $\Theta:U\too U^*$ un isomorphism  which is invariant, i.e., for any $X\in U$,
\begin{equation}
\label{inva}\Theta\circ\mathrm{L}_X=\mathrm{L}_X^*\circ\Theta.
\end{equation}
 We associate to any endomorphism $A:U\too U$ the tensors $\de(A),\mathcal{O}(A)\in\mathrm{End}(U\otimes U,U)$ given by
\begin{equation}\label{ybs}\de(A)(X,Y)=X.A(Y)-Y.A(X)-A([X,Y])\esp \mathcal{O}(A)(X,Y)=[AX,AY]-(A(AX.Y)-A(AY.X)).
\end{equation}One can see easily that
\begin{equation}\label{oeq}\mathcal{O}(A)=N_A+A\circ\de(A),
\end{equation}where $N_A$ is the Nijenhuis torsion of $A$.
The following proposition gives an useful characterization of quasi $S$-matrices and $S$-matrices on $(U,.)$. The second assertion of this proposition was obtained by Bai (see Corollary 6.8 \cite{bai}).
\begin{pr}\label{prqsss}Let $\mathrm{r}\in U\otimes U$ and $\mathfrak{a}$ its skew-symmetric part. Put $A=\mathrm{r}_\#\circ\Theta$ and $T=\mathfrak{a}_\#\circ\Theta$. Then the following assertions holds:
\begin{enumerate}\item[$(i)$] The tensor $\mathrm{r}$ is a quasi $S$-matrix of
$(U,.)$ if and only if $T$ is $\mathrm{L}_U$-invariant and $\mathcal{O}(A)$ is $\mathrm{L}_U^*\otimes\mathrm{L}_U^*\otimes\ad$-invariant.
\item[$(ii)$] If $\mathrm{r}$ is symmetric then it is a  $S$-matrix of
$(U,.)$ if and only if  $\mathcal{O}(A)=0$.
\item[$(iii)$] If $\mathrm{r}$ is symmetric and invertible then it is a  $S$-matrix of
$(U,.)$ if and only if  $\de(A^{-1})=0$. 

 \end{enumerate}
\end{pr}

 {\bf Proof.} We have, for any $X,Y\in U$,
 \begin{eqnarray*}\;\De(\mathrm{r})(\Theta X,\Theta Y)&=&\mathrm{r}_\#([\Theta X,\Theta Y])-
 \;[\mathrm{r}_\#(\Theta X),\mathrm{r}_\#(\Theta Y)]
 \\&\stackrel{\eqref{eq110bis}}=&
 A\circ\Theta^{-1}\left(\mathrm{L}_{AY}^t\Theta X \right)
 -A\circ\Theta^{-1}\left(\mathrm{L}_{AX}^t\Theta Y \right)\\&&+2
 A\circ\Theta^{-1}\left(\mathrm{L}(\mathfrak{a})(\Theta X,\Theta Y) \right)-[AX,AY]\\
 &\stackrel{\eqref{inva}}=&-[AX,AY]-A(AY.X)+A(AX.Y)
 +2
  A\circ\Theta^{-1}\left(\mathrm{L}(\mathfrak{a})(\Theta X,\Theta Y) \right)\\
  &=&-\mathcal{O}(A)(X,Y)+2
    A\circ\Theta^{-1}\left(\mathrm{L}(\mathfrak{a})(\Theta X,\Theta Y) \right).
 \end{eqnarray*}
 From this relation and \eqref{inva} one can see easily that $(i)$ and $(ii)$ holds. Now it is easy that if $A$ is invertible then $\mathcal{O}(A)=0$ if and only if $\de(A^{-1})=0$ and $(iii)$ holds.\hfill$\square$\\
 According to a terminology used by Bai \cite{bai}, if $\mathcal{O}(A)=0$ then $A$ is called un $\mathcal{O}$-operator for the Lie algebra $(U,[\;,\;])$ with respect to the representation $\mathrm{L}_U$.\\
 There are two interesting cases:
 \begin{enumerate}\item[$(i)$] The isomorphism $\Theta$ is skew-symmetric. In this case  $(U,.,\om)$ is a symplectic left symmetric algebra where   $\om(X,Y)=\prec\Theta(X),Y\succ$. 
 \item[$(i)$] The isomorphism $\Theta$ is symmetric. In this case  $(U,.,\prs)$ is a  flat pseudo-Riemannian Lie algebra where   $\langle X,Y\rangle=\prec\Theta(X),Y\succ$.
 
 \end{enumerate}

\begin{pr}\label{prprotheta}Let $(U,.)$ be a left symmetric algebra and $\Theta:U\too U^*$ a invariant isomorphism    and $\mathrm{r}\in U\otimes U$. Put $A=\mathrm{r}_\#\circ\Theta$ and denote by $\circ$ the product on $U$ pull-back  by $\Theta$ of the product on
$U^*$ given by \eqref{eq110}. Then the following assertions holds.
\begin{enumerate}\item If $\Theta$ is skew-symmetric then, for any $X,Y\in U$,
\begin{equation}\label{product} X\circ Y=[AX,Y]+A(Y. X)+Q(X,Y), \end{equation}where $Q:U\times U\too U$ is defined by
\[ \prec\al,Q(X,Y)\succ=-\om(\de(A^s-A^a)(\Theta^{-1}(\al),Y),X),\quad\forall\al\in U^*, \] $A^s$ and $A^a$ are respectively the symmetric part and the skew-symmetric part of $A$ with respect to the 2-form $\om$ associated to $\Theta$.
\item If $\Theta$ is symmetric then, for any $X,Y\in U$, \begin{equation}\label{product1} X\circ Y=Y.AX+AX.Y-A(Y. X)+P(X,Y), \end{equation}where $P:U\times U\too U$ is defined by
\[ \prec\al,P(X,Y)\succ=\langle\de(A^s-A^a)(\Theta^{-1}(\al),Y),X\rangle,\quad\forall\al\in U^*, \]
$A^s$ and $A^a$ are respectively the symmetric part and the skew-symmetric part of $A$ with respect the 2-form $\prs$ associated to $\Theta$.

\item If $\Theta$ is skew-symmetric and $\mathrm{r}$ is a quasi $S$-matrix then $(U,\circ,\om)$ is a symplectic left symmetric algebra if and only if $\de(A^a)=0$.

\end{enumerate}

\end{pr}

{\bf Proof.}\begin{enumerate}\item Suppose that $\Theta$ is skew-symmetric and define $\om$ by $\om(X,Y)=\prec\Theta(X),Y\succ$. So, for any $\al,\be\in U^*$, 
\[ \mathrm{r}(\al,\be)=-\om(A\circ\Theta^{}(\al),\be). \]
We have
\begin{eqnarray*}
\prec\al,X\circ Y\succ&=&\prec\al,\Theta^{-1}(\Theta(X).\Theta(Y))\succ\\
&=&-\prec\Theta(X).\Theta(Y),\Theta^{-1}(\al)\succ\\
&\stackrel{\eqref{eq110}}=&-\mathrm{r}(\mathrm{L}_{\Theta^{-1}(\al)}^t\Theta(X),\Theta(Y))-
\mathrm{r}(\Theta(X),\ad_{\Theta^{-1}(\al)}^t\Theta(Y))\\
&\stackrel{\eqref{inva}}=&\mathrm{r}(\Theta(\Theta^{-1}(\al). X),\Theta(Y) )+\om(AX,\Theta^{-1}\left(\ad_{\Theta^{-1}(\al)}^t\Theta(Y) \right))\\
&=&\om(Y,A(\Theta^{-1}(\al). X))-\om(Y,[\Theta^{-1}(\al),AX])\\
&=&\om((A^s-A^a)Y,\Theta^{-1}(\al). X)-\om(Y,[\Theta^{-1}(\al),AX])\\
&=&-\om(\Theta^{-1}(\al).(A^s-A^a)Y,X)-\om(Y,[\Theta^{-1}(\al),AX])\\
&=&-\om(\de((A^s-A^a)(\Theta^{-1}(\al),Y),X)-\om(Y.(A^s-A^a)(\Theta^{-1}(\al)),X)\\&&-
\om((A^s-A^a)([\Theta^{-1}(\al),Y],X)-\om(Y,[\Theta^{-1}(\al),AX])\\
&=&-\om(\de((A^s-A^a)(\Theta^{-1}(\al),Y),X)+\om(\Theta^{-1}(\al),A(Y.X))-
\om([\Theta^{-1}(\al),Y],AX)-\om(Y,[\Theta^{-1}(\al),AX])\\
&\stackrel{(a)}=&-\om(\de((A^s-A^a)(\Theta^{-1}(\al),Y),X)+
\prec\al,A(Y.X))+[AX,Y]\succ.
\end{eqnarray*}In $(a)$ we have used the fact that $\om$ is 2-cocycle with respect to the Lie bracket.
\item Suppose that $\Theta$ is symmetric and define $\prs$ by $\langle X,Y\rangle=\prec\Theta(X),Y\succ$. So, for any $\al,\be\in U^*$, 
\[ \mathrm{r}(\al,\be)=\langle A\circ\Theta^{}(\al),\be\rangle. \]
We have, for any $\al\in U^*$ and $X,Y\in U$,
\begin{eqnarray*}
\prec\al,X\circ Y\succ&=&\prec\al,\Theta^{-1}(\Theta(X).\Theta(Y))\succ\\
&=&\prec\Theta(X).\Theta(Y),\Theta^{-1}(\al)\succ\\
&\stackrel{\eqref{eq110}}=&\mathrm{r}(\mathrm{L}_{\Theta^{-1}(\al)}^t\Theta(X),\Theta(Y))+
\mathrm{r}(\Theta(X),\ad_{\Theta^{-1}(\al)}^t\Theta(Y))\\
&\stackrel{\eqref{inva}}=&-\mathrm{r}(\Theta(\Theta^{-1}(\al). X),\Theta(Y) )+\langle AX,\Theta^{-1}\left(\ad_{\Theta^{-1}(\al)}^t\Theta(Y) \right)\rangle\\
&=&-\langle Y,A(\Theta^{-1}(\al). X)\rangle+\langle Y,[\Theta^{-1}(\al),AX]\rangle\\
&=&-\langle(A^s-A^a)Y,\Theta^{-1}(\al). X\rangle+\langle Y,[\Theta^{-1}(\al),AX]\rangle\\
&=&\langle\Theta^{-1}(\al).(A^s-A^a)Y,X\rangle+\langle Y,[\Theta^{-1}(\al),AX]\rangle\\
&=&\langle\de((A^s-A^a)(\Theta^{-1}(\al),Y),X\rangle+\langle Y.(A^s-A^a)(\Theta^{-1}(\al)),X\rangle+
\langle(A^s-A^a)([\Theta^{-1}(\al),Y],X\rangle+\langle Y,[\Theta^{-1}(\al),AX]\rangle\\
&=&\langle\de((A^s-A^a)(\Theta^{-1}(\al),Y),X\rangle-\langle\Theta^{-1}(\al),A(Y.X)\rangle+
\langle[\Theta^{-1}(\al),Y],AX\rangle+\langle Y,[\Theta^{-1}(\al),AX]\rangle\\
&=&\langle\de((A^s-A^a)(\Theta^{-1}(\al),Y),X\rangle-\al(A(Y.X))-
\langle Y.\Theta^{-1}(\al),AX\rangle-\langle Y,AX.\Theta^{-1}(\al)\rangle\\
&=&\langle\de((A^s-A^a)(\Theta^{-1}(\al),Y),X\rangle-\prec\al,A(Y.X)\succ+\prec\al,Y.AX+AX.Y\succ.
\end{eqnarray*}
\item Suppose that $\Theta$ is skew-symmetric and $\mathrm{r}$ is a quasi $S$-matrix. We have, for any $X,Y,Z\in U$,
\begin{eqnarray*}
\om(X\circ Y,Z)+\om(Y,X\circ Z)&=&\om([AX,Y],Z)+\om(A(Y.X),Z)+\om(Q(X,Y),Z)\\
&&+\om(Y,[AX,Z])+\om(Y,A(Z. X))+\om(Y,Q(X,Z))\\
&=&-\om(X,(A^s-A^a)[Y,Z]-Y. (A^s-A^a)Z+Z. (A^s-A^a)Y)-\Theta(Z)(Q(X,Y))+\Theta(Y)(Q(X,Z))\\
&=&-\om(X,\de((A^s-A^a))(Y,Z))+\om(\de((A^s-A^a))(Y,Z),X)-\om(\de((A^s-A^a))(Z,Y),X)\\
&=&-3\om(X,\de(A^s-A^a))(Y,Z)).
\end{eqnarray*}To conclude, remark that since $\mathrm{r}$ is a quasi $S$-matrix then its skew-symmetric part is $\mathrm{L}_U$-invariant and hence $\de(A^s)=0$.
\hfill$\square$ 

\end{enumerate}

The proof of the two following theorems is similar to the one of Theorem \ref{theosymp}. The second part of Theorem \ref{theohyper} is based on the third part of Proposition \ref{prprotheta} and \eqref{oeq}.

\begin{theo}\label{theohyper}Let $(U,.,\om)$ be a symplectic left symmetric algebra and $A$ a   endomorphism of $U$. We denote by $A^s$ and $A^a$, respectively, the symmetric and the skew-symmetric part of $A$ with respect to $\om$.  The following assertions hold.
\begin{enumerate}\item If $\mathcal{O}(A)$ is $\mathrm{L}_U^*\otimes\mathrm{L}_U^*\otimes\ad$-invariant and $A^s$ is $\mathrm{L}_U$-invariant
 then:
 \begin{enumerate}\item[$(i)$] the product on $U$ given by \eqref{product} is left symmetric,\item[$(ii)$] $(T(U),[\;,\;]^A,K_A,J_A)$ is a complex product structure and
 $(T(U),[\;,\;]^A,\prs_A,K_A)$  is a para-K\"ahler Lie algebra.\end{enumerate}
 
  \item  If $A^s$ is $\mathrm{L}_U$-invariant, $\de(A^a)=0$ and $N_A$ is 
$\mathrm{L}_U^*\otimes\mathrm{L}_U^*\otimes\ad$-invariant  then
$(T(U),[\;,\;]^A,\prs_A,K_A,J_A)$  is a hyper-para-K\"ahler Lie algebra.
\end{enumerate}On what above, $[\;,\;]^A,\prs_A,K_A, J_A$ are given by
\begin{eqnarray*} [(X,Y),(Z,T)]^A&=&([X,Z]+\mathcal{O}(A)(T,Y),X.T-Z.Y),
\quad J_A(X,Y)=(-Y+AX-A^2Y,X-AY),\\
\langle(X,Y),(Z,T)\rangle_A&=&\om(T,X)+\om(Y,Z)
+2\om(A^aY,T),\quad K_A(X,Y)=(X-2AY,-Y). \end{eqnarray*}

\end{theo}

\begin{theo}\label{theohyper2}Let $(\G,\prs)$ be a flat pseudo-Riemannian Lie algebra and $A$ an   endomorphism of $\G$.  We denote by $A^s$ and $A^a$, respectively, the symmetric and the skew-symmetric part of $A$ with respect to $\prs$.
If $\mathcal{O}(A)$ is $\mathrm{L}_\G^*\otimes\mathrm{L}_\G^*\otimes\ad$-invariant and $A^a$ is $\mathrm{L}_\G$-invariant then the product on $\G$ given by \eqref{product1} is left symmetric. Moreover, $(T(\G),[\;,\;]^A,K_A,J_A)$ is a complex product structure and
 $(T(\G),[\;,\;]^A,\prs_A,K_A)$  is a para-K\"ahler Lie algebra, where
 \begin{eqnarray*} [(X,Y),(Z,T)]^A&=&([X,Z]+\mathcal{O}(A)(T,Y),X.T-Z.Y),
 \quad J_A(X,Y)=(-Y+AX-A^2Y,X-AY),\\
 \langle(X,Y),(Z,T)\rangle_A&=&\langle T,X\rangle+\langle Y,Z\rangle
 +2\langle A^sY,T\rangle,\quad K_A(X,Y)=(X-2AY,-Y). \end{eqnarray*}
 The dot here is the Levi-Civita product and $\mathrm{L}_\G$ is its associated representation.
\end{theo}

\begin{remark}\label{rem4}As in Remark \ref{rem3} we have the following more general result. Let $(U,.)$ be a left symmetric algebra and $A$ an endomorphism of $U$ such that $\mathcal{O}(A)$ is $\mathrm{L}_{U}^*\otimes\mathrm{L}_{U}^*\otimes\ad$-invariant. Then, one can check that the bracket $[\;,\;]^A$ on $T(U)$ is a Lie bracket and hence $L^A:\G\times\G\times\G\too\G$ given by
\[ L^A(X,Y,Z)=\mathcal{O}(A)(X,Y).Z \]is a Lie triple system.

\end{remark}

\section{ Four dimensional hyper-para-K\"ahler Lie algebras}\label{section6}

In this section, we determine all 4-dimensional hyper-para-K\"ahler Lie algebras, up an
isomorphism. To do that we need first to determine two dimensional symplectic left
symmetric algebras and the couples of compatible such algebras. Four dimensional
hyper-para-K\"ahler Lie algebras were classified in  \cite{andrada} by a huge
computation, we adopt here a new method which reduces significantly the calculations.\\

Let $(U,.,\om)$ be a symplectic left symmetric algebra. We have
\begin{eqnarray}
 \label{e1}\mathrm{R}_{u.v}&=&\mathrm{R}_{v}\circ\mathrm{R}_{u}+[\mathrm{L}_{u},\mathrm{R}
_
{v}],\\
\label{e2}\mathrm{L}_{[u,v]}&=&[\mathrm{L}_{u},\mathrm{L}_{v}],\\
\label{e3}\mathrm{L}_{u}&+&\mathrm{L}_u^\mathrm{a}=0,
\end{eqnarray}
where $\mathrm{L}_u^\mathrm{a}$ is the adjoint of $\mathrm{L}_{u}$ with respect to $\om$.
Put
\begin{eqnarray*}
 U.U&=&\mathrm{span}\left\{u.v,u,v\in U\right\},\\
D(U.U)&=&\mathrm{span}\left\{u.v-v.u,u,v\in U\right\},\\
S(U.U)&=&\mathrm{span}\left\{u.v+v.u,u,v\in U\right\}.
\end{eqnarray*}
We have clearly
\begin{equation}\label{eq51}U.U=D(U.U)+S(U.U)\esp (U.U)^\perp=\left\{u\in
U,\mathrm{R}_u=0\right\}.\end{equation}
The sign $\perp$ designs the orthogonal with respect to $\om$.

\begin{pr}\label{pr52}
 Let $(U,.,\om)$ be an abelian symplectic left symmetric algebra. Then
$U.U\subset(U.U)^\perp$.
\end{pr}
{\bf Proof.} We have, for any $u\in U$, $\mathrm{R}_u=\mathrm{L}_u$ and then we get from
\eqref{e1}-\eqref{e2} that, for any $u,v\in U$,
$$\mathrm{R}_{u.v}=\mathrm{R}_u\circ\mathrm{R}_v=\mathrm{R}_v\circ\mathrm{R}_u.$$
Moreover, $\mathrm{R}_u^\mathrm{a}=-\mathrm{R}_u$ so $\mathrm{R}_{u.v}=0$ and the proposition
follows from \eqref{eq51}.
\hfill$\square$
 \begin{pr}\label{pr53}
 Let $(U,.,\om)$ be a 2-dimensional non trivial abelian symplectic left symmetric algebra.
Then
there exists a basis $\{e_1,e_2\}$ of $U$ such that
$$\om=e_1^*\wedge e_2^*,\, \mathrm{R}_{e_1}=\mathrm{L}_{e_1}=0\esp e_2.e_2=ae_1,\quad
a\not=0.$$
\end{pr}
{\bf Proof.} We get from Proposition \ref{pr52} that
$U.U=(U.U)^\perp=\mathrm{span}\{e_1\}$. Choose
$e_2$ such that $\om(e_1,e_2)=1$ and the proposition follows.\hfill $\square$

\begin{pr}\label{pr54}
 Let $(U,.,\om)$ be a 2-dimensional non abelian symplectic left symmetric algebra. Then
there exists a basis $\{e_1,e_2\}$ of $U$ such that
$$\om=e_1^*\wedge e_2^*,\, e_1.e_1=0, e_2.e_2=ae_2\esp e_1.e_2=-e_2.e_1=ae_1,\quad
a\not=0.$$
\end{pr}
{\bf Proof.} We have necessary $\dim D(U.U)=1$. We distinguish two cases:
\begin{enumerate}
 \item {\bf First case}: $\dim U.U=1$. In this case $U.U=D(U.U)=\mathrm{span}\{e_1\}$. If
$S(U.U)=\{0\}$ then we can choose $e_2$ such that $\om(e_1,e_2)=1$. Since
$e_1.e_1,e_2.e_2\in
S(U.U)$
then $e_1.e_1=e_2.e_2=0$. Moreover, since $U.U=(U.U)^\perp$ then $\mathrm{R}_{e_1}=0$. Now
$e_1.e_2\in S(U.U)$ and then $e_1.e_2=0$. Thus the product is trivial.\\
If $S(U.U)\not=\{0\}$ then $U.U=D(U.U)=S(U.U)=\mathrm{span}\{e_1\}$. Choose $e_2$ such
that
$\om(e_1,e_2)=1$. We have $\mathrm{R}_{e_1}=0$, $e_2.e_2=ae_1$ and $e_1.e_2=be_1$. The
relation 
$\mathrm{L}_{[e_1,e_2]}e_2=[\mathrm{L}_{e_1},\mathrm{L}_{e_2}](e_2)$ implies $b=0$ and
hence $[e_1,e_2]=0$ which is 
impossible. In conclusion this case is impossible.
\item {\bf Second case}: $\dim U.U=2$. In this case $U.U=D(U.U)\oplus S(U.U)$. Choose a
basis
$\{e_1,e_2\}$ of $U$ such that $e_1\in D(U.U)$, $e_2\in S(U.U)$ and $\om(e_1,e_2)=1$.
Since
$e_1\in
D(U.U)^\perp$ and $e_2\in S(U.U)^\perp$, we get
\begin{eqnarray*}
 \om(e_1.e_2+e_2.e_1,e_2)&=&0,\\
\om(e_1.e_2+e_2.e_1,e_1)&=&2\om(e_1.e_2,e_1),\\
&=&-2\om(e_2,e_1.e_1)=0.
\end{eqnarray*}Thus
$e_1.e_2=-e_2.e_1$ and
hence $[e_1,e_2]=2e_1.e_2$. So $e_1.e_2=-e_2.e_1=ae_1$. On the other hand,
$e_1.e_1,e_2.e_2\in
S(U.U)$ so
 $e_1.e_1=be_2$ and $e_2.e_2=ce_2$. Now, the relation
$\mathrm{L}_{e_2}^\mathrm{a}=-\mathrm{L}_{e_2}$ implies
$c=a$ and the relation $\mathrm{L}_{[e_1,e_2]}=[\mathrm{L}_{e_1},\mathrm{L}_{e_2}]$
implies $b=0$
and the proposition follows.\hfill$\square$\\
\end{enumerate}

\begin{remark}\label{rem51} From Propositions \ref{pr53} and \ref{pr54}, one can deduce that
if $(U,.,\om)$ is an abelian symplectic left symmetric algebra then $D(U.U)=0$ and
$U.U=S(U.U)$ is an $\om$-isotropic one dimensional vector space. However, if $(U,.,\om)$
is a non abelian symplectic left symmetric algebra then $U=U.U=D(U.U)\oplus S(U.U)$ where 
$D(U.U)$ and $S(U.U)$ are one dimensional $\om$-isotropic vector spaces. This remark will
play a crucial role in the proof of Theorem \ref{theo51}.
 
\end{remark}

Recall that two symplectic left symmetric structures $(U,\star,\om)$ and 
$(U,\circ,\om)$ are called compatible if $K^{\star,\circ}$ satisfies  the second assertion of
Proposition \ref{pr21}. It is obvious that if $(U,\star,\om)$ is a symplectic left
symmetric algebra then it is compatible with $(U,\al\star,\om)$ for any $\al\in\K$. We
call this case trivially compatible.

\begin{theo}\label{theo51}
 Let $(U,\star,\om)$ and 
$(U,\circ,\om)$ be two  symplectic left symmetric structures over a two dimensional
vector space $U$. Then $(U,\star,\om)$ and 
$(U,\circ,\om)$ are non trivially compatible if and only if
one of the following holds:
\begin{enumerate}
 \item There exists a basis $\{e_1,e_2\}$  such that
\begin{eqnarray*}
 {\mathrm{L}}_{e_1}^\star&=&\left(\begin{array}{cc}0&a\\0&0\end{array}\right),\quad
{\mathrm{L}}_{e_2}^\star=\left(\begin{array}{cc}-a&-b\\0&a\end{array}\right),\;
{{\mathrm{L}}}_{e_1}^\circ=0\esp
{{\mathrm{L}}}_{e_2}^\circ=\left(\begin{array}{cc}0&b\\0&0\end{array}\right),
\end{eqnarray*}
$\om=e_1^*\wedge e_2^*$ with $a\not=0$ and $b\not=0$.
\item There exists a basis $\{e_1,e_2\}$  such that
$$\mathrm{L}_{e_1}^\star=\left(\begin{array}{cc}0&a\\0&0\end{array}\right),\;
\mathrm{L}_{e_2}^\star=\left(\begin{array}{cc}-a&b\\0&a\end{array}\right),\;
\mathrm{L}_{e_1}^\circ=\left(\begin{array}{cc}0&c\\0&0\end{array}\right),\;
\mathrm{L}_{e_2}^\circ=\left(\begin{array}{cc}-c&-b\\0&c\end{array}\right).$$
$\om=e_1^*\wedge e_2^*$ with $a\not=0$,  $b\not=0$ and $c\not=0$.

\end{enumerate}

\end{theo}
{\bf Proof.} The proof is based on and adequate use of the fact that the sum of two
compatible symplectic left symmetric structures is left symmetric (see Proposition
\ref{prcompatible}) and the use of Propositions \ref{pr53}-\ref{pr54} and Remark \ref{rem51}. \\

First one can check that  if $\star$ and $\circ$ have one of the form above then they are symplectic and compatible.
Suppose that $(U,\star,\om)$ and 
$(U,\circ,\om)$ are non trivially compatible. We distinguish three cases.
\begin{enumerate}
 \item Both $\star$ and $\circ$ are abelian. Then
$\star+\circ$ defines an abelian symplectic left symmetric algebra structure on $U$
and hence, by virtue of Proposition \ref{pr53},
there exists $a\not=0$ and a basis $\{e_1,e_2\}$ of $U$ such that $\om=e_1^*\wedge e_2^*$ 
$$e_1\star e_1+e_1\circ e_1=e_1\star e_2+e_1\circ e_2=0,\;
e_2\star e_2+e_2\circ e_2=ae_1.\eqno(*)$$ Moreover, $(U\star U)^\perp=U\star U$ and 
$(U\circ U)^\perp=U\circ U$.\\ Suppose that $e_1\circ e_1\not=0$. Then, from $(*)$ above,
we get $U\star U=U\circ U=\mathrm{span}\{e_1\}$,
  and by  \eqref{eq51} $\mathrm{L}_{e_1}^\star=
\mathrm{L}_{e_1}^\circ=0$ which is in contradiction with $e_1\circ e_1\not=0$. So
$e_1\circ e_1=e_1\star e_1=0$. A same argument shows that 
$e_1\circ e_2=e_1\star e_2=0$ and hence $\mathrm{L}_{e_1}^\star=
\mathrm{L}_{e_1}^\circ=0$ which implies by virtue of \eqref{eq51} $U\star U=U\circ
U=\mathrm{span}\{e_1\}$. There exists hence $b\not=0$ and $c\not=0$ such that
$e_2\star e_2=b e_1$ and  $e_2\circ e_2=c e_1$. Finally, $\circ=\frac{b}c\star$ and this
case is not possible.
 \item The product $\star$ is not abelian and $\circ$ is abelian. Then
$\star+\circ$ defines a non abelian symplectic left symmetric algebra structure on $U$
and hence, by virtue of Proposition \ref{pr54},
there exists $a\not=0$ and a basis $\{e_1,e_2\}$ of $U$ such that $\om=e_1^*\wedge e_2^*$
and
$$e_1\star e_1+e_1\circ e_1=0,\; e_2\star e_2+e_2\circ e_2=ae_2,\;
e_1\star e_2+e_1\circ e_2=-e_2\star e_1-e_2\circ e_1=ae_1.\eqno(**)
$$ Moreover, $U=D(U\star U)\oplus S(U\star U)$ and $U\circ U=S(U\circ U)=(U\circ
U)^\perp$.
By adding the  two last relations in $(**)$, we get
$$e_1\star e_2+e_2\star e_1=2 e_1\circ e_2.$$So $e_1\circ e_2\in S(U\star U)$. If
$e_1\circ e_2\not=0$ it spans $U\circ U$ and hence from the second relation in $(**)$ we
deduce that $e_2\in S(U\star U)$ and hence $U\circ U=\mathrm{span}\{e_2\}$. So by
\eqref{eq51}
$\mathrm{L}_{e_2}^\circ=0$ which contradicts $e_1\circ e_2\not=0$. Thus 
$e_1\circ e_2=0$.\\
Suppose now that $e_1\circ e_1\not=0$. We deduce from the second relation in $(**)$ that 
$U\circ U=\mathrm{span}\{e_2\}$ and  then
$\mathrm{L}_{e_2}^\circ=0$. We deduce that
$$e_1\star e_1=-e_1\circ e_1=b e_2,\; e_1\star e_2=-e_2\star e_1=ae_1.$$
From the relation $\mathrm{L}_{[e_1,e_2]}^\star e_1=[\mathrm{L}_{e_1}^\star,
\mathrm{L}_{e_2}^\star]e_1$ we deduce that $b=0$ and $\circ=0$ so we must have 
$e_1\circ e_1=0$.\\ So
we have shown then that $\mathrm{L}_{e_1}^\circ=0$ and hence $U\circ
U=\mathrm{span}\{e_1\}$, $e_2\circ e_2=be_1$. We deduce that
$$e_1\star e_1=0,\;e_2\star e_2=ae_2-be_1,\; e_1\star e_2=-e_2\star e_1=ae_1.$$
We get that $\star$ and $\circ$ satisfy the first form in the theorem.
\item Both $\star$  and $\circ$ are non abelian. According to Proposition
\ref{prcompatible},
$\star+\circ$ and $\star-\circ$ define two symplectic left symmetric algebra
structures on $U$ and one of them must be non abelian. So we can suppose that
$\star+\circ$ is non abelian by replacing $\circ$ by $-\circ$ if it is necessary.
By virtue of Proposition \ref{pr54},
there exists $a\not=0$ and a basis $\{e_1,e_2\}$ of $U$ such that $\om=e_1^*\wedge e_2^*$
and
$$e_1\star e_1+e_1\circ e_1=0,\; e_2\star e_2+e_2\circ e_2=ae_2,\;
e_1\star e_2+e_1\circ e_2=-e_2\star e_1-e_2\circ e_1=ae_1.\eqno(***)
$$ Moreover,  $U=D(U\star U)\oplus S(U\star U)=D(U\circ U)\oplus S(U\circ
U)$. Put
$$v=e_1\star e_2+e_2\star e_1=-e_1\circ e_2-e_2\circ e_1.$$
If $v\not=0$ then it spans $S(U\star U)$ and $S(U\circ U)$ so, from $(***)$ above,
we get $S(U\star U)=S(U\circ U)=\mathrm{span}
\{e_2\}$. So
$$\mathrm{L}_{e_1}^\star=\left(\begin{array}{cc}0&c\\b&0\end{array}\right),\;
\mathrm{L}_{e_2}^\star=\left(\begin{array}{cc}-c&0\\d&c\end{array}\right),\;
\mathrm{L}_{e_1}^\circ=\left(\begin{array}{cc}0&a-c\\-b&0\end{array}\right),\;
\mathrm{L}_{e_2}^\circ=\left(\begin{array}{cc}c-a&0\\-d&a-c\end{array}\right).$$
The relations
$$\mathrm{L}_{[e_1,e_2]}^\star=[\mathrm{L}_{e_1}^\star,\mathrm{L}_{e_2}^\star]
\esp \mathrm{L}_{[e_1,e_2]}^\circ=[\mathrm{L}_{e_1}^\circ,\mathrm{L}_{e_2}^\circ]
\eqno(****)$$are
equivalent to
$d^2=4bc=4b(c-a)$ which is equivalent to $d=b=0$. This implies that
$\circ=\frac{a-c}c\star$. This case is impossible and hence 
$v=0$. \\
If $w=e_1\star e_1=-e_1\circ e_1\not=0$, then it spans $S(U\star U)$ and $S(U\circ U)$ so,
from $(***)$,
we get $S(U\star U)=S(U\circ U)=\mathrm{span}
\{e_2\}$. So
$$\mathrm{L}_{e_1}^\star=\left(\begin{array}{cc}0&c\\b&0\end{array}\right),\;
\mathrm{L}_{e_2}^\star=\left(\begin{array}{cc}-c&0\\0&c\end{array}\right),\;
\mathrm{L}_{e_1}^\circ=\left(\begin{array}{cc}0&a-c\\-b&0\end{array}\right),\;
\mathrm{L}_{e_2}^\circ=\left(\begin{array}{cc}c-a&0\\0&a-c\end{array}\right).$$
In this case, $(****)$ imply  $b=0$ and hence $\circ=\frac{a-c}c\star$. This
case is
impossible and hence $w=0$. To summarize, we have shown that
$$e_1\star e_2+e_2\star e_1=-e_1\circ e_2-e_2\circ e_1=
e_1\star e_1=-e_1\circ e_1=0.$$So
$$\mathrm{L}_{e_1}^\star=\left(\begin{array}{cc}0&c\\0&0\end{array}\right),\;
\mathrm{L}_{e_2}^\star=\left(\begin{array}{cc}-c&d\\0&c\end{array}\right),\;
\mathrm{L}_{e_1}^\circ=\left(\begin{array}{cc}0&a-c\\0&0\end{array}\right),\;
\mathrm{L}_{e_2}^\circ=\left(\begin{array}{cc}c-a&-d\\0&a-c\end{array}\right).$$
We get that $\star$ and $\circ$ satisfy the second form in the theorem.\\
Finally, a direct computation shows that for $\star$ and $\circ$ as in the first form in
the theorem, we have
$$K^{\star,\circ}(e_1,e_1)=K^{\star,\circ}(e_1,e_2)=
K^{\star,\circ}(e_2,e_1)=0\esp K^{\star,\circ}(e_2,e_2)=
-2\left(\begin{array}{cc}0&ab\\0&0\end{array}\right),$$
and for the second form
$$K^{\star,\circ}(e_1,e_1)=K^{\star,\circ}(e_1,e_2)=
K^{\star,\circ}(e_2,e_1)=0\esp K^{\star,\circ}(e_2,e_2)=
2\left(\begin{array}{cc}0&ab+bc\\0&0\end{array}\right).$$
In the two cases
$K^{\star,\circ}$ satisfies  the  second assertion of Proposition \ref{pr21} and the
theorem is proved.\hfill
$\square$
\end{enumerate}

\section{Symplectic associative algebras}\label{section7}
In this section, we study 	an important subclass of the class of symplectic left symmetric algebras. In order to
introduce this subclass  we begin by giving a geometric interpretation of symplectic left
symmetric algebras.\\
Let $(U,.,\om)$ be a symplectic left symmetric algebra. The  product being
Lie-admissible, the bracket $[u,v]=u.v-v.u$ is a Lie bracket on $U$. Moreover, since 
$$\om(u.v,w)+\om(v,u.w)=0$$ for any $u,v,w\in U$, a direct computation gives
$$\om([u,v],w)+\om([v,w],u)+\om([w,u],v)=0,$$ and hence $(U,[\;,\;],\om)$ is a symplectic Lie
algebra. Let $G$ be the simply connected Lie group associated to $(U,[\;,\;])$. For any $u\in
U$ denote by $u^\ell$ the left invariant vector field on $G$ associated to $u$. The formulas
$$\om^\ell(u^\ell,v^\ell)=\om(u,v)\esp \na_{u^\ell}v^\ell=\left(u.v\right)^\ell$$define on $G$
a
left invariant symplectic form and a  flat and torsion free left invariant  connection. Moreover,
$\na$ is symplectic, i.e., $\na\om^\ell=0$. The connection $\na$ is right invariant iff, for
any $u,v,w\in U$,
$$[u^\ell,\na_{v^\ell}w^\ell]=\na_{[u^\ell,v^\ell]}w^\ell+\na_{v^\ell}[u^\ell,w^\ell].$$
A straightforward computation shows that this relation is equivalent to the associativity of
the left symmetric product on $U$. \\
A {\it symplectic associative algebra} is a  symplectic left symmetric algebra which is
associative. We have seen that there is a correspondence between the set of symplectic
associative algebras and symplectic Lie groups endowed with a bi-invariant affine structure
for which the symplectic form is parallel.

In what follows we will give an accurate description of symplectic associative algebras (see Theorems \ref{main1}-\ref{main2}).

 Let $(U,.,\om)$ be an associative
symplectic algebra. Then, for any $u,v\in U$,
$$\mathrm{L}_{uv}=\mathrm{L}_{u}\circ\mathrm{L}_{v}.$$
Since $\mathrm{L}_{u}^\mathrm{a}=-\mathrm{L}_{u}$ for any $u\in U$, we get
\begin{equation}\label{s7eq1}\mathrm{L}_{uv}=\mathrm{L}_{u}\circ
\mathrm{L}_{v}=-\mathrm{L}_{v}
\circ\mathrm{L}_{u}=-\mathrm{L}_{vu}.\end{equation}
\begin{pr}\label{s7pr1}Let $(U,.,\om)$ be an associative
symplectic algebra. Then $U^4=0$ and $\mathcal{J}=U^2+(U^2)^\perp$ is a co-isotropic two-side
ideal of $U$ satisfying  $\mathcal{J}^2=0$. 
 
\end{pr}

{\bf Proof.}  It is obvious that $\mathcal{J}$ is a co-isotropic two-side ideal.  For any $u,v,w\in U$,
\begin{eqnarray*} 
\mathrm{L}_{uv w}&=&\mathrm{L}_{u}\circ\mathrm{L}_{v}\circ\mathrm{L}_{w} \\
&\stackrel{\eqref{s7eq1}}=&\mathrm{L}_{w}\circ\mathrm{L}_{u}
\circ\mathrm{L}_{v}\\
&=&\mathrm{L}_{w}\circ\mathrm{L}_{uv}\\
&\stackrel{\eqref{s7eq1}}=&-\mathrm{L}_{uv}\circ\mathrm{L}_{w}\\
&=&-\mathrm{L}_{uv w},
\end{eqnarray*}
 hence $\mathrm{L}_{uv w}=0$ and $U^4=0$.  Recall that $$(U^2)^\perp=\left\{u\in U, \mathrm{R}_{u}=0 \right\}.$$ So $U^2.(U^2)^\perp=(U^2)^\perp.(U^2)^\perp=0$. On the other hand, for any $u\in (U^2)^\perp$ and any $v,w\in U$,
 $$u.v.w\stackrel{\eqref{s7eq1}}=-v.u.w=0,$$ so $(U^2)^\perp.U^2=0$. In conclusion, $\mathcal{J}^2=0$.
\hfill $\square$\\

According to this proposition, to study  associative symplectic algebras we distinguish two cases depending on the triviality of $U^3$ or not. 

\paragraph{Model of associative symplectic algebras with $U^3=0$}  Let $V$ be a vector space and $(I,\mathfrak{s})$ a symplectic vector space. Let $\mathfrak{m}: V^*\too V\odot V$ and $\mathfrak{n}:I\too V\odot V$ be two linear maps ($V\odot V$ is the space of bilinear symmetric forms on $V^*$). \\
The space $U_1=V\oplus I\oplus V^*$ carries a symplectic form $\om$ for which $I$ and $V\oplus V^*$ are orthogonal, $\om_{|I\times I}=\mathfrak{s}$ and for any $u\in V$, $\al\in V^*$, $\om(\al,u)=-\om(u,\al)=\al(u)$. Define a product on 
$U$ such that $V^*.V^*,I.V^*\subset V$ by
\begin{equation} \label{model1}\prec\ga,\al.\be\succ=\mathfrak{m}(\al)(\be,\ga)\esp \prec\be,i.\al\succ=\mathfrak{n}(i)(\al,\be), \end{equation}
for any $\al,\be,\ga\in V^*$ and $i\in I$ (all the others products vanish).\\
It is easy to see that $(U_1,.,\om)$ is an associative symplectic algebra and $U_1^3=0$. We call such algebra {\it associative symplectic algebra of type one.} Actually, all associative symplectic algebras with $U^3=0$ are of this form.
\begin{theo}\label{main1}Any   associative symplectic algebra with $U^3=0$ is isomorphic to an associative symplectic algebra of type one.
\end{theo}

{\bf Proof.} The condition $U^3=0$ is equivalent to $U^2\subset(U^2)^\perp$. Put $V=U^2$ and choose $I$ a complement of $V$ in $(U^2)^\perp$. The restriction of $\om$ to  $I$ defines a  symplectic form say $\mathfrak{s}$. The orthogonal $I^\perp$ of $I$ is a symplectic space which  contains $V$ as a Lagrangian subspace so we can choose a Lagrangian complement $W$ of $V$ in $I^\perp$. The linear map $W\too V^*$, $u\mapsto\om(u,.)$ is an isomorphism so we can identify $(U,\om)$ to $V\oplus I\oplus V^*$ endowed with the symplectic form described above. Since $V.U=U.V=0$ and for any $u\in I$, $\mathrm{R}_u=0$ we get that the product on $U$ is given by \eqref{model1} which achieves the proof.\hfill$\square$

\paragraph{Model of associative symplectic algebras with $U^3\not=0$}  Let $V=V_0\oplus V_1$  and $I=I_0\oplus I_1$ be two  vector spaces such that
$(I_0,\mathfrak{s_0})$, $(I_1,\mathfrak{s_1})$  are symplectic vector spaces. Put $V^*=V_0^0\oplus V_1^0$ where $V_0^0$ and $V_1^0$ are the annihilators of $V_0$ and $V_1$ respectively.\\
Let $\mathfrak{a}: V_1\too V_0\odot V_0$,  $\mathfrak{b}:I_0\too V_0\odot V_0$,  $\mathfrak{c}:I_1\too V\odot V$ and $\mathfrak{d}:V^*\too V\odot V$   linear maps ($V_0\odot V_0$ (resp. $V_1\odot V_1$) is the space of bilinear symmetric forms on $V_1^0$ (resp. $V_0^0$)). Finally, let $F:V^*\times V_1^0\too I_0$ a bilinear map.\\
The space $U_2=V\oplus I\oplus V^*$ carries a symplectic form $\om$ for which $I$ and $V\oplus V^*$ are orthogonal, $\om_{|I\times I}=\mathfrak{s_0}\oplus\mathfrak{s_1}$ and for any $u\in V$, $\al\in V^*$, $\om(\al,u)=-\om(u,\al)=\al(u)$.\\
Define now a product on $U_2$ satisfying 
$$V_1.V_1^0,I_0.V_1^0,V^*.I_0 \subset V_0,\: I_1.V^*\subset V, V^*.V^0_0\subset V,\: V^*.V^0_1\subset V\oplus I_0$$and given by
\begin{eqnarray*}
\prec\be_1,v_1.\al_1\succ&=&\mathfrak{a}(v_1)(\al_1,\be_1),\:v_1\in V_1,\al_1,\be_1\in V_1^0,\\
\prec\be_1,i_0.\al_1\succ&=&\mathfrak{b}(i_0)(\al_1,\be_1),\:i_0\in I_0,\al_1,\be_1\in V_1^0,\\
\prec\be_1,\al.i_0\succ&=&\mathfrak{s}_0(i_0,F(\al,\be_1)),\:i_0\in I_0,\al\in V^*,\be_1\in V_1^0,\\
\prec\be,i_1.\al\succ&=&\mathfrak{c}(i_1)(\al,\be),\:i_1\in I_1,\al,\be\in V^*,\\
\al.\be_0&=&E(\al,\be_0),\: \al\in V^*,\be_0\in V_0^0,\\
\al.\be_1&=&E(\al,\be_1)+F(\al,\be_1),\: \al\in V^*,\be_1\in V_1^0.
\end{eqnarray*} 
\begin{eqnarray*} 
\prec\ga,E(\al,\be)\succ=\mathfrak{d}(\al)(\be,\ga).\end{eqnarray*}
With this product $U_2$ becomes an algebra for which the multiplication by left is symplectic. Now this product is associative iff, for any $\al\in V^*$, $\be_1,\ga_1,\mu_1\in V_1^0$, $\be_0\in V_0^0$, 
\begin{eqnarray}\label{eqs71}
\mathfrak{a}{(E_1(\al,\be_1))}(\ga_1,\mu_1)+\mathfrak{b}{(F(\al,\be_1))}(\ga_1,\mu_1)=
\mathfrak{s}_0(F(\be_1,\ga_1),F(\al,\mu_1)),\\\label{eqs72}
\mathfrak{a}{(E_1(\al,\be_0))}(\ga_1,\mu_1)=\mathfrak{s}_0(F(\be_0,\ga_1),F(\al,\mu_1)).
\end{eqnarray}
In this case $U_2^3=0$ iff
\[ \mathfrak{s}_0(F(\be_1,\ga_1),F(\al,\mu_1))=\mathfrak{s}_0(F(\be_0,\ga_1),F(\al,\mu_1))=0. \]
When \eqref{eqs71} and \eqref{eqs72} hold and $U_2^3\not=0$ we call $(U_2,.,\om)$ an {\it associative symplectic algebra of type two.}
\begin{theo}\label{main2}Any associative symplectic algebra with $U^3\not=0$ is isomorphic to an associative algebra of type two.
\end{theo}

{\bf Proof.} Since $U^4=0$ and $U^3\subset U^2$ then
$$U^3\subset U^2\subset (U^3)^\perp\esp U^3\subset (U^2)^\perp\subset (U^3)^\perp.$$
Put $V_0=U^3$ and choose a complement $V_1$ of $V_0$ in $V=U^2\cap(U^2)^\perp$. Choose $I_0$ and $I_1$
two subspaces of $U$ such that
$$U^2=V\oplus I_0\esp (U^2)^\perp=V\oplus I_1.$$
We have $I_0\cap I_1=\left\{0\right\}$, $\om(I_0,I_1)=0$, $I_0,I_1$ are symplectic and  $I=I_0\oplus I_1$ is also
symplectic.  Denote by $\mathfrak{s}_0$ and $\mathfrak{s}_1$ the restrictions of $\om$ to $I_0$ and $I_1$ respectively. 
Now $I^\perp$ is symplectic and contains $V$ as a Lagrangian subspace so we can choose a
Lagrangian subspace $W$ of $I^\perp$ complement of $V$. The linear map $W\too V^*$,
$u\mapsto\om(u,.)$ realizes an isomorphism. Finally, we get 
the identification 
$$U=V\oplus(I_1\oplus I_2)\oplus V^*$$with the symplectic form given by
$$\om(V,V)=\om(V,I)=\om(V^*,V^*)=\om(V^*,I)=0\esp \om_{|I\times I}=\om_1\oplus\om_2,$$ and for any $u\in V$, $\al\in V^*$, $\om(\al,u)=-\om(u,\al)=\al(u).$ Denote by $V_0^0$ and $V_1^0$ the annihilator of $V_0$ and $V_1$ respectively.\\
Let us study the product's properties. We have shown in Proposition \ref{s7pr1}
that $\mathcal{J}=V\oplus(I_0\oplus I_1)$ satisfies $\mathcal{J}^2=0$. We have obviously that for any $u\in V_0$, $\mathrm{L}_u=\mathrm{R}_u=0$ and
 from \eqref{eq51}  for any $u\in V\oplus I_1$, $\mathrm{R}_u=0$. Since $U^2.V^*\subset V_0$ and the fact that the symplectic form is invariant we get $U^2.V_0^0=0$. Since $V^*.U^2\subset V_0$ we get from the invariance of the symplectic form that $V^*.V_0^0\subset V$. So we can put
 \begin{eqnarray*}
\prec\be_1,v_1.\al_1\succ&=&\mathfrak{a}(v_1)(\al_1,\be_1),\:v_1\in V_1,\al_1,\be_1\in V_1^0,\\
\prec\be_1,i_0.\al_1\succ&=&\mathfrak{b}(i_0)(\al_1,\be_1),\:i_0\in I_0,\al_1,\be_1\in V_1^0,\\
\prec\be_1,\al.i_0\succ&=&\mathfrak{s}_0(i_0,F(\al,\be_1)),\:i_0\in I_0,\al\in V^*,\be_1\in V_1^0,\\
\prec\be,i_1.\al\succ&=&\mathfrak{c}(i_1)(\al,\be),\:i_1\in I_1,\al,\be\in V^*,\\
\al.\be_0&=&E(\al,\be_0),\: \al\in V^*,\be_0\in V_0^0,\\
\al.\be_1&=&E(\al,\be_1)+F(\al,\be_1),\: \al\in V^*,\be_1\in V_1^0.
\end{eqnarray*} The invariance of symplectic form implies that $\mathfrak{a}(v_1)$, $\mathfrak{b}$, $\mathfrak{c}(i_1)$ and $\mathfrak{d}(\al)$ are symmetric, where 
$\mathfrak{d}(\al)(\be,\ga)=\ga(E(\al,\be))$. The associativity of this product is equivalent to $\al.(\be.\ga)=(\al.\be).\ga$ for any $\al,\be,\ga\in V^*$. This is obviously true when $\ga\in V_0^0$. So the associativity is equivalent to
$$\mu_1(\al.(\be_0.\ga_1))=\mu_1((\al.\be_0).\ga_1)\esp
\prec\mu_1,\al.(\be_1.\ga_1)\succ=\prec\mu_1,(\al.\be_1).\ga_1\succ$$for any
$\al\in V^*$, $\al_1,\be_1,\mu_1\in V_1^0$ and $\be_0\in V_0^0$. Which is equivalent to
\begin{eqnarray*}
\mathfrak{s}_0(F(\be_0,\ga_1),F(\al,\mu_1))&=&
\mathfrak{a}(E_1(\al,\be_0))(\ga_1,\mu_1),\\
\mathfrak{s}_0(F(\be_1,\ga_1),F(\al,\mu_1))&=&
\mathfrak{a}(E_1(\al,\be_1))(\ga_1,\mu_1)
+\mathfrak{b}(F(\al,\be_1))(\ga_1,\mu_1).
\end{eqnarray*}This achieves the proof.\hfill$\square$

\begin{co}Let $(U,.,\om)$ be an associative symplectic algebra. Then
\begin{itemize}\item[(i)]If $\dim U=2$ then $(U,.,\om)$ is isomorphic to an associative symplectic algebra of type one of the form $V\oplus V^*$ with $\dim V=1$.
\item[(ii)]If $\dim U=4$ then $(U,.,\om)$ is isomorphic to an associative symplectic algebra of type one either of the form $V\oplus V^*$ with $\dim V=2$
or $V\oplus I\oplus V^*$ with $\dim V=1$.

\end{itemize}

\end{co}

A six dimensional associative algebra $U^3\not=0$ is isomorphic to $V\oplus I_0\oplus V^*$ with $\dim V=2$ and $V=V_0\oplus V_1$. Choose a basis $(e_0,e_1)$ a basis of $V$ such $e_i\in V_i$ and $(f_1,f_2)$ a basis of $I_0$ such that $\mathfrak{s}_1(f_1,f_2)=1$. The equations \eqref{eqs71} and \eqref{eqs72} are equivalent to 
\begin{eqnarray*}
\mathfrak{a}{(E_1(e_0^*,e_1^*))}(e_1^*,e_1^*)+\mathfrak{b}{(F(e_0^*,e_1^*))}
(e_1^*,e_1^*)=
\mathfrak{s}(F(e_1^*,e_1^*),F(e_0^*,e_1^*)),\\
\mathfrak{a}{(E_1(e_0^*,e_0^*))}(e_1^*,e_1^*)=\mathfrak{s}(F(e_0^*,e_1^*),F(e_0^*,e_1^*)),\\
\mathfrak{a}{(E_1(e_1^*,e_1^*))}(e_1^*,e_1^*)+\mathfrak{b}{(F(e_1^*,e_1^*))}
(e_1^*,e_1^*)=
\mathfrak{s}(F(e_1^*,e_1^*),F(e_1^*,e_1^*)),\\
\mathfrak{a}{(E_1(e_1^*,e_0^*))}(e_1^*,e_1^*)=\mathfrak{s}(F(e_0^*,e_1^*),F(e_1^*,e_1^*)).
\end{eqnarray*}
Put
\begin{eqnarray*}
F(e_0^*,e_1^*)&=&af_1+bf_2\esp
F(e_1^*,e_1^*)=cf_1+df_2.
\end{eqnarray*}
So
\begin{eqnarray*}
\mathfrak{a}{(E_1(e_0^*,e_1^*))}(e_1^*,e_1^*)+a\mathfrak{b}{(f_1)}
(e_1^*,e_1^*)+b\mathfrak{b}{(f_2)}
(e_1^*,e_1^*)=cb-ad,\\
\mathfrak{a}{(E_1(e_0^*,e_0^*))}(e_1^*,e_1^*)=0,\\
\mathfrak{a}{(E_1(e_1^*,e_1^*))}(e_1^*,e_1^*)+c\mathfrak{b}{(f_1)}
(e_1^*,e_1^*)+d\mathfrak{b}{(f_2)}
(e_1^*,e_1^*)=0,\\
\mathfrak{a}{(E_1(e_1^*,e_0^*))}(e_1^*,e_1^*)=ad-cb.
\end{eqnarray*}
Which is equivalent to
\begin{eqnarray*}
a\mathfrak{b}{(f_1)}
(e_1^*,e_1^*)+b\mathfrak{b}{(f_2)}
(e_1^*,e_1^*)&=&-\mathfrak{a}{(E_1(e_1^*,e_0^*))}(e_1^*,e_1^*)-\mathfrak{a}{(E_1(e_0^*,e_1^*))}(e_1^*,e_1^*),\\
c\mathfrak{b}{(f_1)}
(e_1^*,e_1^*)+d\mathfrak{b}{(f_2)}
(e_1^*,e_1^*)&=&-\mathfrak{a}{(E_1(e_1^*,e_1^*))}(e_1^*,e_1^*),\\
\mathfrak{a}{(E_1(e_0^*,e_0^*))}(e_1^*,e_1^*)&=&0,\\
\mathfrak{a}{(E_1(e_1^*,e_0^*))}(e_1^*,e_1^*)&=&ad-cb.
\end{eqnarray*}
Put $\mathfrak{a}{(e_1)}(e_1^*,e_1^*)=\al\not=0$, $E_1(e_0^*,e_0^*)=0$ and
 $\de=ad-cb\not=0$, $E_1(e_i^*,e_j^*)=a_{ij}e_1$. So 
 \[\mathfrak{b}{(f_1)}
 (e_1^*,e_1^*)=-\al\de^{-1}(aa_{11}-c(a_{10}+a_{01}))\esp \mathfrak{b}{(f_2)}
 (e_1^*,e_1^*)=\al\de^{-1}(ba_{11}-d(a_{10}+a_{01})).  \]

\bibliographystyle{elsarticle-num}

\end{document}